\documentclass[reqno,12pt]{amsart}

\usepackage{times}
\usepackage{lmodern}
\usepackage{textcomp}
\usepackage{amsmath,cancel}
\usepackage{amsfonts,enumitem}
\usepackage{amssymb}
\usepackage{color}
\usepackage{mathrsfs}
\usepackage[dvipsnames]{xcolor}
\usepackage{esint}
\usepackage[colorlinks=true]{hyperref}
\usepackage{cleveref}
\newtheorem{thm}{Theorem}[section]
\newtheorem{lem}[thm]{Lemma}
\newtheorem{cor}[thm]{Corollary}
\newtheorem{prop}[thm]{Proposition}

\numberwithin{equation}{section}

\newcommand{\R}{\mathbb{R}}

\newcommand{\N}{\mathbb{N}}

\newcommand{\ve}{\varepsilon}
\newcommand{\rd}{\mathrm{d}}

\newcommand{\dhr}{\mathrel{\lhook\joinrel\relbar\kern-.8ex\joinrel\lhook\joinrel\rightarrow}}

\newcommand{\cgray}[1]{{\color{gray} #1}}

\begin{document}

\title[Chemotaxis model with indirect signal production]{Identifying the critical mass in a quasilinear chemotaxis model with indirect signal production}\thanks {\texttt{orcid: 0000-0003-3091-8085 (PhL), 0000-0002-2017-5049 (ChW)}}

\author{Philippe Lauren\c{c}ot}
\address{Laboratoire de Math\'ematiques (LAMA) UMR~5127, Universit\'e Savoie Mont Blanc, CNRS\\	F--73000 Chamb\'ery, France}
\email{philippe.laurencot@univ-smb.fr}
\author{Christoph Walker}
\address{Leibniz Universit\"at Hannover\\
Institut f\"ur Angewandte Mathematik\\
Welfengarten 1\\
30167 Hannover\\
Germany}
\email{walker@ifam.uni-hannover.de}
\date{\today}

\begin{abstract}
The mass threshold separating boundedness and unboundedness of global solutions to a critical quasilinear parabolic-elliptic chemotaxis model with indirect signal production is identified in the radially symmetric setting in a ball. Specifically, solutions emanating from initial values with mass below a critical mass $M_c>0$, which is characterized in terms of a variant of the Hardy-Littlewood-Sobolev inequality, are global and bounded. In contrast, for any $M>M_c$, there are initial values with mass $M$ such that the corresponding global solutions are unbounded. The proof relies on the variational structure of the system and a detailed study of the stationary states.
\end{abstract}
%
\maketitle

\section{Introduction}\label{sec.1}

Chemotaxis systems describe the dynamics of a population of cells with density $u\ge 0$, undergoing diffusion and chemotactic movement in response to a chemical signal with concentration $v\ge 0$, which is produced directly by the cells or indirectly (by an intermediate species, for instance). For the so-called \textsl{critical mass chemotaxis systems,} a typical feature is the existence of a sharp mass threshold $\mathscr{M}_{th}>0$, which separates two different dynamical behaviors, according to the value of the initial mass $\|u(0)\|_1$ of the cell density. More precisely, when $\|u(0)\|_1$ ranges in $(0,\mathscr{M}_{th})$, the corresponding solution is global and bounded. In contrast, for any $M>\mathscr{M}_{th}$ (possibly except for a discrete set), there are initial values with $\|u(0)\|_{1} = M$ for which the solution becomes unbounded in finite or infinite time. A prototypical example is the celebrated Keller-Segel system
\begin{equation}
\begin{aligned}
	\partial_t u  & = \mathrm{div}\big(\nabla u - u\nabla v\big)\,, \qquad & t>0\,, \ x\in \Omega\,,\\
	\tau \partial_t v & = \Delta v - v + u \,, \qquad & t>0\,, \ x\in\Omega\,,
\end{aligned} \label{lks}
\end{equation}
in space dimension $d=2$, with $\mathscr{M}_{th}=8\pi$ when $\Omega$ is a ball in $\mathbb{R}^2$ (or is $\mathbb{R}^2$) in the radially symmetric setting, and $\mathscr{M}_{th}=4\pi$ when $\Omega$ is a smooth bounded domain of $\mathbb{R}^2$ and $(u,v)$ satisfies homogeneous Neumann boundary conditions. Similar results hold for the parabolic-elliptic Keller-Segel system (corresponding to $\tau=0$ in~\eqref{lks}), and also for the J\"ager-Luckhaus and Smoluchowski-Poisson versions, see \cite{Bi1998, Bi2020, DP2004, GZ1998, Ho2002, Ho2003, HW2001, JL1992, NSY1997} and the references therein. For~\eqref{lks}, this dichotomy is a genuinely two-dimensional phenomenon and cannot be observed otherwise. Indeed, all solutions are global and bounded in one space dimension \cite{OY2001}, while there are initial values $u(0)$ with arbitrary positive mass $\|u(0)\|_1$ leading to finite-time blowup in space dimensions $d\ge 3$ \cite{Wi2013}. Nevertheless, in other chemotaxis systems, featuring either a nonlinear diffusion or a non-constant chemotactic sensitivity, a  mass threshold phenomenon may occur in space dimensions $d\ge 3$. For instance, the existence of such a threshold is reported in \cite{CG2023, Wi2022} for the Keller-Segel chemotaxis system with volume-filling effect
\begin{equation}
\begin{aligned}
	\partial_t u  & = \mathrm{div}\big( D(u)\nabla u - u S(u)\nabla v\big)\,, \qquad & t>0\,, \ x\in \Omega\,,\\
	\tau \partial_t v & = \Delta v - v + u \,, \qquad & t>0\,, \ x\in\Omega\,,
\end{aligned} \label{vfks}
\end{equation}
with diffusion coefficient $D\ge 0$ and chemotactic sensitivity $S\ge 0$, when the ratio $D/S$ satisfies, roughly speaking,
\begin{equation}
	\frac{D}{S}(z) \sim C z^{m_c-1} \;\;\text{  as }\;\; z\to\infty\,, \qquad m_c := \frac{2(d-1)}{d}\,. \label{ras}
\end{equation} 
In that case, a threshold mass $\mathscr{M}_{th}^b$ is defined in an abstract way as
\begin{subequations}\label{McA}
\begin{equation}
	\mathscr{M}_{th}^b := \sup\left\{ M>0\ :\ \begin{array}{l}
		\text{the solution is global and bounded}\\
		\text{whenever } u(0)\in X \;\text{ with }\; \|u(0)\|_1\le M
	\end{array} \right\}\,, \label{McAa}
\end{equation}
where $X$ is an appropriate functional subspace of $L_1(\Omega)$, and is shown to be positive and finite in \cite{CG2023, Wi2022}. Since one also expects unbounded solutions, we may define in an analogous way another threshold mass $\mathscr{M}_{th}^u$ as 
\begin{equation}
	\mathscr{M}_{th}^u := \inf\left\{ M>0\ :\ \begin{array}{l}
		\text{for any } \tilde{M}>M \text{ (possibly except for a discrete}\\
        \text{set), there is } u(0)\in X \text{ with } \|u(0)\|_1 = \tilde{M} \text{ such}\\
        \text{that the solution is unbounded }
	\end{array} \right\}\,. \label{McAb}
\end{equation}
\end{subequations}
Clearly, $\mathscr{M}_{th}^b\le \mathscr{M}_{th}^u$, and we shall say that there is a critical mass threshold whenever they coincide, putting $\mathscr{M}_{th} := \mathscr{M}_{th}^b =\mathscr{M}_{th}^u$ in that case. As discussed above, $\mathscr{M}_{th}=4\pi$ for~\eqref{lks} for a general smooth bounded domain~$\Omega$ in $\mathbb{R}^2$ and is thus explicitly determined. Explicitly identifying $\mathscr{M}_{th}^b$ and $\mathscr{M}_{th}^u$ remains largely an open problem and, to the best of our knowledge, is only achieved in a handful of cases for the specific choice
\begin{equation}
	D(z)=m_c z^{m_c-1}\,, \quad S(z)=1\,, \qquad z\ge 0\,, \label{qks}
\end{equation} 
with $m_c$ defined in~\eqref{ras}. More precisely, introducing the Newton potential 
\begin{equation}
	E_d(x) :=c_d |x|^{2-d}\,, \quad x\in\mathbb{R}^d\setminus\{0\}\,, \qquad c_d:=\frac{1}{(d-2)\sigma_d}\,, \label{Np}
\end{equation}
where $\sigma_d := 2\pi^{d/2}/\Gamma(d/2)$ is the surface area of the unit sphere~$\mathbb{S}^{d-1}$, it is shown in \cite{BCL2009, ST2009} that, for the critical quasilinear Smoluchowski-Poisson equation in $\mathbb{R}^d$, $d\ge 3$,
\begin{equation}
	\partial_t u = \mathrm{div}\big(\nabla u^{m_c} - u \nabla (E_d*u)\big)\,, \qquad t>0\,,\ x\in\mathbb{R}^d \,,\label{qsp}
\end{equation} 
there holds $\mathscr{M}_{th}^b = \mathscr{M}_{th}^u$ and the critical mass $\mathscr{M}_{th}$ equals $M_c$ defined in~\eqref{Mc} below, see also \cite{OW2019}. In addition, all solutions to~\eqref{vfks} with $D$ and $S$ given by~\eqref{qks} and $d\ge 3$ are global and bounded as soon as $\|u(0)\|_1\in (0,M_c)$ when $\Omega=\mathbb{R}^d$ \cite{BL2013, Mi2017}. When $\Omega$ is a ball in $\mathbb{R}^d$ (or is $\mathbb{R}^d$) and $M>M_c$, finite-time blowup is shown in the radially symmetric setting for some initial values $u(0)$ with $\|u(0)\|_1=M$, but only for $d\in\{3,4\}$ \cite{LM2017}. Consequently, $\mathscr{M}_{th}=M_c$ is known to hold in that case only when $\Omega=\mathbb{R}^d$ and $d\in\{3,4\}$, even though it is expected to be true for all $d\ge 3$. We emphasize that, in contrast to the two-dimensional case, the threshold mass $M_c$ is not given explicitly. Nevertheless, it is characterized by its connection to the optimal constant in a functional inequality, see~\eqref{Mc}, or equivalently by the solution to an elliptic equation, see~\eqref{b0}.

The existence of a threshold mass is also investigated for the Keller-Segel chemotaxis system with indirect signal production
\begin{equation}
\begin{aligned}
	\partial_t u  & = \mathrm{div}\big( D(u)\nabla u - u S(u)\nabla v\big)\,, \qquad & t>0\,, \ x\in \Omega\,,\\
	\tau \partial_t v & = \Delta v - \delta v + w - (1-\delta) \langle w\rangle \,, \qquad & t>0\,, \ x\in \Omega\,,\\
	\partial_t w & = u-w\,, \qquad & t>0\,, \ x\in \Omega\,,
\end{aligned} \label{klip}
\end{equation}
where $(\tau,\delta)\in \{(0,0),(0,1),(1,1)\}$, supplemented with $(1-\delta) \langle v \rangle = 0$. When $(\tau,\delta)=(0,0)$ and $\Omega$ is a ball in $\mathbb{R}^2$, it is shown in \cite{TW2017} for the radially symmetric setting that $\mathscr{M}_{th}=8\pi$ when $D=S\equiv1$, while the positivity and finiteness of $\mathscr{M}_{th}^b$, along with the finiteness of $\mathscr{M}_{th}^u$, is obtained in \cite{FLT2023} when $(D,S)$ are given by~\eqref{qks}. Subsequently, we proved in \cite{LW26} that these values coincide and identified them as equal to $M_c$. The threshold masses $\mathscr{M}_{th}^b$ and $\mathscr{M}_{th}^b$ are also positive and finite when 
\begin{equation*} 
    D\equiv 1 \;\;\text{ and }\;\; S(z) = (1+z)^{1-m_c}\,, \quad z\ge 0\,,
\end{equation*}
but for a smaller class of initial values, see \cite{DLY2026}. We point out that the approach used in \cite{DLY2026, FLT2023, TW2017} relies on the special form of the elliptic equation satisfied by $v$ when $(\tau,\delta)=(0,0)$, which allows one to reduce the system for $(u,v)$ to a single parabolic equation for the cumulative distribution of $u$ for which a comparison principle is available. However, this technique appears to yield an identification of $\big(\mathscr{M}_{th}^b,\mathscr{M}_{th}^b\big)$ only  when $D=S\equiv 1$ \cite{TW2017}. Deriving the identity $\mathscr{M}_{th}=M_c$ in \cite{LW26} requires a different approach of a more variational flavor, which is more in the spirit of the one developed in \cite{Ho2002, HW2001} for the Keller-Segel chemotaxis system~\eqref{lks}. Moreover, the reduction to a single parabolic equation does not seem to extend to the cases $(\tau,\delta)\in\{(0,1),(1,1)\}$, and the variational approach described above is used in \cite{La2019} to show that $\mathscr{M}_{th}=4\pi$ when $(\tau,\delta)\in \{(0,1),(1,1)\}$, $D=S\equiv 1$, and $\Omega$ is a smooth bounded domain of $\mathbb{R}^2$.

\medskip

The main purpose of this work is to prove the existence of a critical mass threshold $\mathscr{M}_{th}$ for~\eqref{klip} with $(\tau,\delta) = (0,1)$ in the radially symmetric setting, when $\Omega$ is a ball in $\mathbb{R}^d$, $d\ge 3$, and $(D,S)$ are given by~\eqref{qks}, recalling that $m_c=2(d-1)/d$. Furthermore, we identify $\mathscr{M}_{th}$ as equal to $M_c$ defined in~\eqref{Mc}. Thus, from now on, we assume that
\begin{equation}
	d\ge 3\,, \quad m=m_c = \frac{2(d-1)}{d}\,, \quad \Omega = \mathbb{B}_R\,, \label{mc}
\end{equation} 
where $\mathbb{B}_R$ denotes the ball in $\R^d$ of radius $R>0$ centered at the origin, and focus on the system
\begin{subequations}\label{E}
\begin{align}
	&\partial_t u  = \mathrm{div}\big(m(1+u)^{m-1}\nabla u-u\nabla v\big)\,, \qquad &(t,x)\in (0,\infty)\times\mathbb{B}_R\,, \label{E1}\\
	&-\Delta v + v = w\,, \qquad &(t,x)\in (0,\infty)\times\mathbb{B}_R\,, \label{E2}\\
	&\partial_t w=u-w\,, \qquad &(t,x)\in (0,\infty)\times\mathbb{B}_R\,, \label{E3} \\
	&\nabla u\cdot \mathbf{n}  = \nabla v\cdot \mathbf{n} = 0\,, \qquad &(t,x)\in (0,\infty)\times\partial\mathbb{B}_R\,, \label{E4} \\
	&(u,w)(0)  = (u^0,w^0)\,, \qquad &x\in \mathbb{B}_R\,. \label{E5}
\end{align}
\end{subequations}
We  present our main results regarding the existence and identification of $\mathscr{M}_{th}$ for~\eqref{E} in Theorem~\ref{THM1} below.

\medskip

\paragraph{\textbf{Notation}} For $p\in (1,\infty)$, we denote by $\Delta_N$ the Laplacian with domain
\begin{equation*}
	W_{p,N}^2(\mathbb{B}_R) := \{ f\in W_p^2(\mathbb{B}_R)\, :\, \nabla f\cdot \mathbf{n} = 0 \;\text{ on }\; \partial\mathbb{B}_R\}.
\end{equation*}
We also introduce the isomorphism
\begin{equation}\label{opK}
	\mathcal{K}:L_p(\mathbb{B}_R)\rightarrow  W_{p,N}^2(\mathbb{B}_R)\,,\quad z\mapsto (1-\Delta_N)^{-1} z\,.
\end{equation} 
Note that $\mathcal{K}[z]\ge 0$ for $z\ge 0$ and that
\begin{equation}\label{int}
\int_{\mathbb{B}_R} \mathcal{K}[z](x)\,\rd x= \int_{\mathbb{B}_R} z(x)\,\rd x\,,\quad z\in L_p(\mathbb{B}_R)\,.
\end{equation}
The subset of radially symmetric functions of a set $X$ of functions defined on $\mathbb{B}_R$ is denoted by $X_{rad}$. We slightly abuse notation by writing $f(x)$, $x\in \bar{\mathbb{B}}_R$, and $f(r)$, $r = |x|\in [0,R]$, simultaneously for radially symmetric functions $f$ defined on the closed ball~$\bar{\mathbb{B}}_R$; i.e. we identify radially symmetric functions with their profiles. The identity
\begin{equation}
	\int_0^R r^{d-1} f(r)\,\rd r = \frac{\|f\|_1}{\sigma_d} = \frac{R^d}{d} \langle f\rangle\,, \qquad f\in L_{1,rad}^+(\mathbb{B}_R)\,, \label{eqi}
\end{equation}
is often used in the following.

\section{Main results}\label{sec.mr}

We begin with the global well-posedness of~\eqref{E} established in \cite{LW26+}. 

\begin{thm}\label{T1:Ex}
Let $p\in (2d,\infty)$ and consider an initial value
\begin{equation}\label{a1}
    (u^0,w^0)\in W_{p}^{1,+}(\mathbb{B}_R)  \times L_\infty^+(\mathbb{B}_R) \,.
\end{equation}
Then there is a unique couple $(u,w)$ of functions satisfying
\begin{subequations}
\begin{align}
&u \in  C^1\big((0,\infty),L_p(\mathbb{B}_R)\big) \cap  C\big((0,\infty),W_{p,N}^{2}(\mathbb{B}_R)\big)\cap C\big([0,\infty),W_{p}^{1}(\mathbb{B}_R) \big)\,,\label{regularityA}\\
& w\in C^1 \big([0,\infty),L_\infty(\mathbb{B}_R)\big)\,,\qquad v := \mathcal{K}[w] \in C^1\big([0,\infty),W_{p,N}^{2}(\mathbb{B}_R)\big)\,,\label{regularityB}
\end{align}
such that $(u,v,w)$ is a global non-negative strong solution to~\eqref{E}.
\end{subequations}
In fact, the mapping 
\begin{equation*}
	\boldsymbol{\Psi}:(u^0,w^0)\mapsto (u,w)
\end{equation*} 
defines a global semiflow on $W_{p}^{1,+}(\mathbb{B}_R) \times L_\infty^+(\mathbb{B}_R)$.  
Moreover, 
\begin{itemize}
    \item  for $t\ge 0$, it holds that
    \begin{equation}
        \|u(t)\|_1 = \|u^0\|_1\,, \quad \|v(t)\|_1 = \|u^0\|_1\,, \quad \|w(t)\|_1 \le \max\big\{ \|w^0\|_1 , \|u^0\|_1 \big\}\,; \label{L1}
    \end{equation}
	\item  if $u\in L_\infty((0,\infty),L_m(\mathbb{B}_R))$, then the orbit $(u,w)([0,\infty))$ is relatively compact in \mbox{$W_p^1(\mathbb{B}_R) \times L_\infty(\mathbb{B}_R)$};
	\item  if $(u^0,w^0)\in W_{p,rad}^{1,+}(\mathbb{B}_R)  \times L_{\infty,rad}^+(\mathbb{B}_R)$, then $(u,v,w)(t)$ is radially symmetric for all $t\ge 0$.
\end{itemize}
\end{thm}

The local well-posedness is based on the theory developed for quasilinear parabolic initial boundary value problems and outlined in \cite{Am1988, Am1989, Am1990, Am1993}, see the references therein, while global existence is achieved by estimates derived along the lines of \cite{FLT2023, TW2017}.

\medskip 

From now on, given an initial value $(u^0,w^0)$ satisfying~\eqref{a1}, we shall refer to 
\begin{equation*}
    (u,w)=\boldsymbol{\Psi}\big(u^0,w^0\big)
\end{equation*} 
as the (corresponding) solution to~\eqref{E}, with the implicit definition 
\begin{equation*}
    v := \mathcal{K}[w]\,.
\end{equation*}
A salient feature of~\eqref{E} is the availability of a Lyapunov functional, a property already observed in \cite{La2019} for $m=1$ and shared with the classical chemotaxis system, see \cite[Theorem~2]{Ho2001} for instance. Specifically, defining $\phi:(0,\infty)\rightarrow (0,\infty)$  by
\begin{equation}\label{phi}
	\phi''(z)=m\frac{(1+z)^{m-1}}{z}\,,\quad z>0\,,\qquad \phi(1)=\phi'(1)=0\,,
\end{equation}
the functional 
\begin{equation*}
	\mathcal{L}: L_m^+(\mathbb{B}_R)\times  L_p^+(\mathbb{B}_R)\rightarrow \mathbb{R}\,,
\end{equation*}
given by
\begin{equation}\label{FL}
	\mathcal{L}(u,w):=\int_{\mathbb{B}_R}\big(\phi(u)- u\, \mathcal{K}[w]\big)\,\rd x+\frac{1}{2}\int_{\mathbb{B}_R} |\nabla\mathcal{K}[w]|^2\,\rd x +\frac{1}{2}\int_{\mathbb{B}_R}\vert\mathcal{K}[w]\vert^2\,\rd x\,,
\end{equation}
is well-defined for all $(u,w) \in  L_m^+(\mathbb{B}_R)\times  L_p^+(\mathbb{B}_R)$ since $p> 2d$ and is a strict Lyapunov functional for~\eqref{E}, as reported in the next result. 

\begin{prop}\label{PP2}
Let $p\in (2d,\infty)$ and assume that the initial value $(u^0,w^0)$ satisfies~\eqref{a1}. Then the corresponding solution $(u,w)$ to~\eqref{E} satisfies
\begin{equation*}
\mathcal{L}(u,w)(t)+\int_0^t \mathcal{D}(u,w)(s)\,\rd s\le \mathcal{L}\big(u^0,w^0\big)\,,\quad t\ge 0\,,
\end{equation*}
where
\begin{equation*}
	\mathcal{D}(u,w):=\int_{\mathbb{B}_R} u\big\vert \nabla\big(\phi'(u)-\mathcal{K}[w]\big)\big\vert^2\,\rd x+\|\nabla \mathcal{K}[u-w]\|_2^2+\|\mathcal{K}[u-w]\|_2^2\ge 0\,.
\end{equation*}
In particular, $\mathcal{L}$ is a strict Lyapunov functional for~\eqref{E}. 
\end{prop}

We are now in a position to state our main result, for which we recall from Theorem~\ref{T1:Ex} that, for $ M>0$, the set 
\begin{equation}
    \mathcal{I}_{M,rad}(\mathbb{B}_R) := \left\{ 
    (u,w)\in W_{p,rad}^{1,+}(\mathbb{B}_R) \times L_{\infty,rad}^+(\mathbb{B}_R)\,:\, \|u\|_1=M 
	\right\} \label{IS}
\end{equation}
is positively invariant for the semiflow $\boldsymbol{\Psi}$.
 
\begin{thm}\label{THM1}
Let $X=W_{p,rad}^{1,+}(\mathbb{B}_R) \times L_{\infty,rad}^+(\mathbb{B}_R)$ and assume that $\mathscr{M}_{th}^b$ and $\mathscr{M}_{th}^u$ are defined in~\eqref{McA}. Then $\mathscr{M}_{th}^b= \mathscr{M}_{th}^u = M_c>0$ with $M_c$ defined in~\eqref{Mc} below. More precisely:
\begin{itemize}
\item [\textbf{(I)}] For any initial value $(u^0,w^0)\in \mathcal{I}_{M,rad}(\mathbb{B}_R)$ with $M\in (0,M_c)$, the global solution $(u,w)$ to~\eqref{E} is bounded.
\item [\textbf{(II)}] Let $M>M_c$. There is $\mu_{M}^s\in\mathbb{R}$ such that the set
\begin{equation*}
    \mathcal{F}_{M} := \left\{ (u^0,w^0)\in \mathcal{I}_{M,rad}(\mathbb{B}_R)\ :\ \mathcal{L}(u^0,w^0)<\mu_{M}^s \right\}
\end{equation*}
is non-empty and, for any initial value $(u^0,w^0)\in \mathcal{F}_{M}$, the global solution $(u,w)$ to~\eqref{E} blows up in infinite time.
\end{itemize}
\end{thm}

As mentioned above, Theorem~\ref{THM1} is already known for the J\"ager-Luckhaus version of~\eqref{E}, which corresponds to the choice $(\tau,\delta)=(0,0)$ and $(D,S)$ given by~\eqref{qks} in~\eqref{klip}, see \cite{LW26}, as well as \cite{TW2017} for the two-dimensional case. As already mentioned, when $v$ is given by
\begin{equation}
    - \Delta v = w - \langle w\rangle\,, \quad \langle v \rangle = 0\,, \qquad (t,x)\in (0,\infty)\times\mathbb{B}_R\,, \label{F}
\end{equation}
instead of~\eqref{E2}, the specific structure of~\eqref{E1} combined with~\eqref{F} allows one to reduce the system for the pair $(u,v)$ to a single parabolic equation for the cumulative distribution function associated with $u$ and use comparison arguments. The analysis performed in \cite{DLY2026, FLT2023, TW2017} relies on this property. Such an approach does not seem to be available when $u$ and $v$ are coupled through~\eqref{E1}-\eqref{E2}, so that a different argument is needed. We instead take advantage of the variational structure of~\eqref{E} as in \cite{Ho2002, HW2001, LW26} and proceed as follows: Given $\bar{M}\ge M>0$, we set
\begin{equation}
    \mathcal{I}_{M,\bar{M},rad}(\mathbb{B}_R) := \left\{ 
    (u,w)\in \mathcal{I}_{M,rad}(\mathbb{B}_R)\,:\, \|w\|_1 = \bar{M} \right\} \label{ISb}
\end{equation}
and denote by $\mathcal{S}_{M,rad}$ the set of stationary solutions to~\eqref{E} in $\mathcal{I}_{M,\bar{M},rad}(\mathbb{B}_R)$, observing that any $(u,w)\in \mathcal{S}_{M,rad}$ satisfies $\|u\|_1=\|w\|_1=M$. We then prove that
\begin{subequations}\label{infl}
\begin{equation}
    \inf_{\mathcal{S}_{M,rad}} \mathcal{L} = \inf_{\mathcal{I}_{M,\bar{M},rad}(\mathbb{B}_R)} \mathcal{L} > -\infty \label{infl1}
\end{equation}
for $M\in (0,M_c)$, while
\begin{equation}
    \mu_M^s := \inf_{\mathcal{S}_{M,rad}} \mathcal{L} > \inf_{\mathcal{I}_{M,\bar{M},rad}(\mathbb{B}_R)} \mathcal{L} = -\infty \label{infl2}
\end{equation}
\end{subequations}
for $M>M_c$, for sufficiently large $\bar{M}$. When $M\in (0,M_c)$, condition~\eqref{infl1} yields a suitable coercivity property for $\mathcal{L}$, which ensures a time-independent bound on $u$ in $L_m(\mathbb{B}_R)$. Theorem~\ref{T1:Ex}  then leads to Theorem~\ref{THM1}~\textbf{(I)}. For $M>M_c$, we follow the strategy developed in \cite{Ho2002, HW2001} and subsequently employed in \cite{CG2023, HW2005, IOY2013, LW26, Wi2010, Wi2022}. It relies on the following consequence of the compactness of the orbits established in Theorem~\ref{T1:Ex} and the Lyapunov functional: every solution to~\eqref{E} which is bounded in $L_\infty(\mathbb{B}_R,\mathbb{R}^2)$ stabilizes to stationary states by LaSalle's invariance principle. Consequently, boundedness is excluded by~\eqref{infl2} when $(u^0,w^0)\in \mathcal{I}_{M,\bar{M},rad}(\mathbb{B}_R)$ satisfies $\mathcal{L}(u^0,w^0)<\mu_M^s$. Thus, the main building blocks of our analysis are properties~\eqref{infl}. First, the radial symmetry of the solution $v$ to~\eqref{E2} implies that it has constant boundary values $V:=v(x)$, $x\in\partial\mathbb{B}_R$, so that $v-V\in H_0^1(\mathbb{B}_R)$. Its trivial extension then belongs to~$W_2^1(\mathbb{R}^d)$. We can therefore compare $\mathcal{L}$ with a similar functional defined on $\mathbb{R}^d$, which is bounded from below according to the variant~\eqref{VHLS} of the Hardy-Littlewood-Sobolev inequality when $M\in (0,M_c)$. The proof of~\eqref{infl2} is more involved and combines scaling arguments with refined elliptic estimates in the spirit of \cite{LW26,LS1994,WY2003} to rule out the existence of an unbounded sequence in $\mathcal{S}_{M,rad}$ provided $M\ne M_c$.   

\medskip

\paragraph{\textbf{Outline.}} We devote the next section to the precise characterization of the critical mass $M_c$ and related results. Proposition~\ref{PP2} and Theorem~\ref{THM1}~\textbf{(I)} are proved in Section~\ref{sec.3}, where also the infimum of $\mathcal{L}$ over $\mathcal{I}_{M,\bar{M},rad}(\mathbb{B}_R)$ is studied, establishing the corresponding statements in~\eqref{infl}. The core of this work is Section~\ref{sec.4}, where the finiteness of $\mu_M^s$ is shown for every $M\ne M_c$. To this end, we study the limiting behavior of sequences $(u_l,w_l)_{l\ge 1}$ in $\mathcal{S}_{M,rad}$ such that $\mathcal{L}(u_l,w_l)\to -\infty$ as $l\to\infty$. A first step is to show that $\|u_l\|_\infty\to \infty$ and that $(u_l)_{l\ge 1}$ behaves as the Dirac mass $M \delta_0$ as $l\to\infty$. At this stage, information on the concentration profile and the behavior of the mass outside the origin are lacking. Thus, a refined study is required to identify both, relying in particular on scaling arguments and elliptic estimates in the spirit of \cite{LS1994, WY2003}.

\medskip

Let us finally mention that the proof of our main result, Theorem~\ref{THM1}, has a similar structure as that of \cite[Theorem~2.3]{LW26} with some arguments being identical. For the reader's convenience, however, we have retained some of those arguments in the proof herein.

\section{The critical mass \texorpdfstring{$M_c$}{}}\label{sec.cmr}

In order to characterize the critical mass $M_c$, we recall that a variant of the  Hardy–Littlewood–Sobolev inequality \cite[Lemma~3.2]{BCL2009} asserts that the constant
\begin{equation}\label{VHLS}
	C_*:=\frac{1}{c_d}\sup\left\{ \dfrac{\displaystyle\int_{\mathbb{R}^d}h(x) (E_d\ast h)(x)\,\rd x}{\|h\|_{L_1(\mathbb{R}^d)}^{2/d}\,\|h\|_{L_m(\mathbb{R}^d)}^m}\ :\, h\in L_1(\mathbb{R}^d)\cap L_m(\mathbb{R}^d)\,,\, h\not= 0 \right\}
\end{equation}
is finite and positive. We set
\begin{equation}\label{Mc}
	M_c:=\left[\frac{2}{(m-1)c_d C_*}\right]^{d/2}\in (0,\infty)\,,
\end{equation}
which is also the critical mass for the quasilinear Patlak–Keller–Segel model with critical nonlinear diffusion~\eqref{qsp} \cite{BCL2009}, as already mentioned. An alternative characterization of $M_c$ stems from the existence of minimizers to~\eqref{VHLS}. In fact, 
\begin{subequations}\label{b0}
\begin{equation}\label{b2x}
	M_c = \int_{\mathbb{B}_1}\zeta^{1/(m-1)}(x)\,\rd x\,,
\end{equation}
 where $\zeta$ is the unique radially symmetric classical solution  to 
\begin{equation}\label{b1}
	-\Delta \zeta = \frac{m-1}{m} \zeta^{1/(m-1)}\quad \text{ in }\ \mathbb{B}_1\,,\qquad \zeta=0\quad \text{ on }\ \partial\mathbb{B}_1\,.
\end{equation}
\end{subequations}
We collect additional properties of $\zeta$ in the next result.

\begin{prop}\label{propB1}
The potential $E_d*\zeta^{1/(m-1)}$ of $\zeta^{1/(m-1)}$ is given by
\begin{equation}
    \big(E_d*\zeta^{1/(m-1)}\big)(x) = \frac{m}{m-1} \zeta(x) + \frac{2-m}{m-1} \frac{\mathscr{I}_\zeta}{M_c}\,, \quad x\in\mathbb{B}_1\,, \label{eq.B2}
\end{equation}
and
\begin{equation}
   \big(E_d*\zeta^{1/(m-1)}\big)(x) = \frac{M_c}{(d-2)\sigma_d} \frac{1}{|x|^{d-2}}\,, \quad x\in\mathbb{R}^d\setminus\mathbb{B}_1\,, \label{eq.B3} 
\end{equation}
where
\begin{equation}
    \mathscr{I}_\zeta := \int_{\mathbb{B}_1} \zeta^{m/(m-1)}(x)\,\rd x\,. \label{eq.B4}
\end{equation}
Moreover,
\begin{equation}
    \mathscr{I}_\zeta = \frac{m-1}{2-m}\frac{M_c^2}{(d-2)\sigma_d} = \frac{m}{m-1} \int_{\mathbb{B}_1} |\nabla\zeta(x)|^2\,\rd x\,. \label{eq.B5}
\end{equation}
\end{prop}

We refer to \cite[Proposition~B.1]{LW26} for a proof.

\section{Proof of Theorem~\ref{THM1}~\textbf{(I)}}\label{sec.3}

In this section, we focus on the functional $\mathcal{L}$ defined in~\eqref{FL}. We first prove that it is indeed a strict Lyapunov functional for~\eqref{E}, as stated in Proposition~\ref{PP2}, and then investigate its boundedness or unboundedness from below. We begin with the continuity of $\mathcal{L}$ on $L_m^+(\Omega)\times  L_p^+(\Omega)$  for $p>2d$.

\begin{lem}\label{lem.contL} 
Let $p\in (2d,\infty)$. Then $\mathcal{L}$ belongs to $C\big(L_m^+(\Omega)\times  L_p^+(\Omega)\big)$. 
\end{lem}

\begin{proof}
This follows exactly as in \cite[Lemma~2.1]{LW26} from the continuous embedding of $W_p^2(\Omega)$ into $W_2^1(\Omega)\cap L_{m/(m-1)}(\Omega)$.
\end{proof}

\begin{proof}[Proof of Proposition~\ref{PP2}]
For simplicity, we only provide a formal proof, assuming implicitly that $u$ is positive on $[0,\infty)\times\bar{\Omega}$, in order not to be bothered by the singularity of $\phi'(z)$ as $z\to 0$. We refer to the proof of \cite[Proposition~2.2]{LW26} for a detailed argument. We derive from~\eqref{E} and \eqref{phi} that 
\begin{align*}
	\frac{\rd}{\rd t}\mathcal{L}(u,w)&= \int_\Omega\big(\phi'( u)- \mathcal{K}[w]\big)\partial_t u\,\rd x-\int_\Omega u\partial_t \mathcal{K}[w]\,\rd x\\
	&\quad+\int_\Omega \nabla \mathcal{K}[w]\cdot\nabla\partial_t \mathcal{K}[w]\,\rd x+\int_\Omega  \mathcal{K}[w] \partial_t \mathcal{K}[w]\,\rd x\\
	&= -\int_\Omega u \big\vert\nabla\big(\phi'( u)- \mathcal{K}[w]\big)\big\vert^2\,\rd x -\int_\Omega (u-w)\partial_t \mathcal{K}[w]\,\rd x\\
	&= -\int_\Omega u \big\vert\nabla\big(\phi'( u)- \mathcal{K}[w]\big)\big\vert^2\,\rd x -\int_\Omega (u-w)\mathcal{K}[u-w]\,\rd x\,.
\end{align*}
Since~\eqref{opK} and Gauß' theorem ensure
\begin{equation*}
	\int_\Omega (w-u)\mathcal{K}[u-w]\,\rd x=\|\mathcal{K}[u-w]\|_2^2+\|\nabla\mathcal{K}[u-w]\|_2^2\,,
\end{equation*}
integration with respect to time yields the claim.
\end{proof}

Next, as in \cite[Lemma~4.2]{LW26}, we take advantage of the radially symmetric setting to control the boundary value of $v=\mathcal{K}[w]$.

\begin{lem}\label{LL1}
Let $w\in L_{1,rad}^+(\mathbb{B}_R)$. Then $v=\mathcal{K}[w]$ is radially symmetric and  
\begin{equation*}
	0\le v(R) \le c_0\|w||_1 \quad \text{with}\quad  c_0:=\frac{2d+R^2}{2\sigma_d R^d}\,.
\end{equation*}
\end{lem}

\begin{proof}
The radial symmetry of $v$ is obvious due to~\eqref{opK} and the rotational invariance of the Laplace operator.  Since equation~\eqref{opK} reads in radial variable
\begin{equation*}
	\partial_r\left(r^{d-1}\partial_rv(r)\right)=r^{d-1}\left( v(r)- w(r) \right)\,,\quad 0<r<R\,,
\end{equation*}
we deduce first 
\begin{equation}\label{hiii}
	\partial_r v(r) =r^{1-d}\int_0^r s^{d-1}\big(v(s)-w(s)\big)\,\rd s\,,\quad 0<r\le R\,,
\end{equation}
so that, invoking the Neumann boundary condition for $v$,
\begin{equation}\label{hii}
	\int_0^R s^{d-1}v(s) \,\rd s=\int_0^R s^{d-1}w(s) \,\rd s\,.
\end{equation}
Integrating then~\eqref{hiii} and using~\eqref{hii}, it follows that
\begin{align*}
	v(R)-v(r)&=-\int_r^R \rho^{1-d}\int_{\rho}^R s^{d-1}\big(v(s)-w(s)\big)\,\rd s \rd \rho\\
	&=\frac{1}{d-2}\int_r^R \big(s-s^{d-1}r^{2-d}\big)\big(v(s)-w(s)\big)\,\rd s \,.
\end{align*}
Therefore, multiplying both sides by $r^{d-1}$ and integrating with respect to $r\in (0,R)$ entail
\begin{align*}
	v(R) & = \frac{d}{R^d}\int_0^R r^{d-1}v(r)\,\rd r-\frac{1}{2R^d}\int_0^R s^{d+1}\big(v(s)-w(s)\big)\,\rd s\,.
\end{align*}
Recalling that $v\ge 0$ and invoking~\eqref{eqi} and again~\eqref{hii} yield the assertion.
\end{proof}

Next, we prepare a coercivity property of $\mathcal{L}$ when $M\in (0,M_c]$.
 
\begin{lem}\label{LLL2}
Let $(u,w)\in  L_{m,rad}^+(\mathbb{B}_R)\times  L_{p,rad}^+(\mathbb{B}_R)$ be such that $M=\|u\|_1\le M_c$ and set $v= \mathcal{K}[w]$. Then
\begin{equation*}
	\mathcal{L}(u,w)\ge \frac{1}{2}\left(1-\left(\frac{M}{M_c}\right)^{2/d}\right) \|\nabla v\|_2^2  +\frac{1}{2}\|v||_2^2- c_0M\| w\|_1 -\frac{m M}{m-1}\,,
\end{equation*}
where $c_0>0$ is the constant from Lemma~\ref{LL1}.
\end{lem}

\begin{proof}
The proof is basically the same as \cite[Lemma~4.3]{LW26}. Indeed, setting
 $V:= v(R)$ we infer as therein from~\eqref{VHLS},~\eqref{Mc}, and the properties of the Newton potential~\eqref{Np} that
\begin{align}
	\int_{\mathbb{B}_R} u v\,\rd x & \le \frac{1}{2}\, \left(\frac{M}{M_c}\right)^{2/d}\,\|\nabla v\|_{2}^2+ \frac{1}{m-1}\|u\|_{m}^{m}\, +VM\,,\label{r101}
\end{align}
and using~\eqref{eq.A2}, we obtain
\begin{align*}
	\mathcal{L}(u,w)&=\int_{\mathbb{B}_R} \phi(u)\,\rd x-\int_{\mathbb{B}_R} u v\,\rd x+\frac{1}{2}\|\nabla v\|_{2}^2+\frac{1}{2}\| v\|_{2}^2\\
	&\ge \frac{1}{2} \left(1-\left(\frac{M}{M_c}\right)^{2/d}\right)\,\|\nabla v\|_{2}^2 +\frac{1}{2}\| v\|_{2}^2- VM  -\frac{mM}{m-1} \,.
\end{align*}
The assertion now follows from the estimate on $V=v(R)$ established in Lemma~\ref{LL1}.
\end{proof}

Given $\bar M\ge M>0$ we now compute the infimum of $\mathcal{L}$ on the set $\mathcal{I}_{M,\bar M,rad}(\mathbb{B}_R)$ defined in~\eqref{ISb}.

\begin{prop}\label{prop1}
Defining
\begin{equation}
    \mu_{M,\bar{M}} := \inf_{(u,w)\in \mathcal{I}_{M,\bar{M},rad}(\mathbb{B}_R)} \mathcal{L}(u,w) \in [-\infty,\infty)\,, \quad \bar{M}\ge M>0\,, \label{mu}
\end{equation}
there is $\bar{M}_0>M$ depending only on $d$, $R$ and $M$ such that
\begin{equation}
\begin{split}
    \mu_{M,\bar{M}} & = -\infty \;\;\text{ for }\;\; M>M_c \;\text{ and }\; \bar{M}\ge \bar{M}_0\,,\\ 
    \mu_{M,\bar{M}} & > -\infty \;\;\text{ for }\;\; M\in (0,M_c] \;\text{ and }\;\bar{M}\ge M\,.
\end{split}\label{x1}
\end{equation}
More precisely, in the second case, the lower bound
\begin{equation}
\begin{split}
    \mathcal{L}(u,w) \ge \frac{1}{2} \left( 1 - \left( \frac{M}{M_c} \right)^{2/d} \right) \|\nabla\mathcal{K}[w]\|_2^2 +\frac{1}{2} \|\mathcal{K}[w]\|_2^2 - c_0 M\bar{M}- \frac{mM}{m-1}\,. \label{lbL}
\end{split}
\end{equation}
holds true for all $(u,w)\in \mathcal{I}_{M,\bar{M},rad}(\mathbb{B}_R)$, recalling that $c_0>0$ is the constant defined in Lemma~\ref{LL1}.
\end{prop}

\begin{proof}
For $M\in (0,M_c]$ and $\bar{M}\ge M$, Proposition~\ref{prop1} is a consequence of Lemma~\ref{LLL2} and the embedding of $W_p^1(\mathbb{B}_R)$ into $L_m(\mathbb{B}_R)$.

Consider the case $\bar{M}\ge M>M_c$. As in  \cite[Proposition~4.4]{LW26} we adapt the analysis of~\cite[Proposition~3.4]{BCL2009} to our setting including  boundary conditions. Let $\kappa:= M/M_c>1$ and define, for $\lambda>\max\{1,1/R\}$ and $\zeta$ given in~\eqref{b1},
\begin{subequations}\label{fl}
\begin{align}
    u_{\lambda}(x) & := \kappa \lambda^d \zeta^{1/(m-1)}(\lambda x)\,, \qquad x\in \mathbb{B}_R\,, \label{ul}\\
    v_{\lambda}(x) & := \kappa \lambda^{d-2} \big(E_d*\zeta^{1/(m-1)}\big)(\lambda x) + \frac{M}{2\sigma_d R^d} |x|^2+\frac{M}{|\mathbb{B}_R|}\,, \qquad x\in \mathbb{B}_R\,, \label{Vl} \\
    w_{\lambda}(x) & := -\Delta v_{\lambda}(x)+v_{\lambda}(x)\,, \qquad x\in\mathbb{B}_R\,. \label{wl}
\end{align}
\end{subequations}
Note first from~\eqref{Vl} and \eqref{wl} that, for $x\in \mathbb{B}_R$,
\begin{align*}
	w_{\lambda}(x) &= \kappa \lambda^{d}\zeta^{1/(m-1)}(\lambda x) + \kappa \lambda^{d-2} \big(E_d*\zeta^{1/(m-1)}\big)(\lambda x)+ \frac{M}{2\sigma_d R^d} |x|^2\ge 0\,.
\end{align*}
We show that there exists $\bar M_0=\bar M_0(d,R,M)>0$ such that $(u_{\lambda},w_{\lambda})$ belongs to $\mathcal{I}_{M,\bar{M},rad}(\mathbb{B}_R)$ whenever $\bar M\ge \bar M_0$. Clearly, $u_{\lambda}$ is a radially symmetric function in $W_p^{1,+}(\mathbb{B}_R)$ with compact support in $\bar{\mathbb{B}}_{1/\lambda}\subset \mathbb{B}_R$ owing to $\lambda R>1$ and $\mathrm{supp}\, \zeta = \bar{\mathbb{B}}_1$. Moreover,~\eqref{b2x} readily entails that
\begin{equation}
	\|u_\lambda\|_1   = \frac{M}{M_c} \int_{\mathbb{B}_{1}} \zeta^{1/(m-1)}(y)\,\rd y = M\,. \label{ub01}
\end{equation}
Next, $v_{\lambda}$ is also a radially symmetric function and belongs to $W_\infty^2(\mathbb{B}_R)$ so that $w_\lambda\in L_{\infty,rad}^+(\mathbb{B}_R)$. Furthermore, according to~\eqref{eq.B3} and $\lambda>1$,
\begin{align}
	\|w_\lambda\|_1&=\|v_\lambda\|_1=\frac{\kappa}{\lambda^2}\int_{\mathbb{B}_{\lambda R}}\big(E_d*\zeta^{1/(m-1)}\big)(y)\,\rd y + \frac{M}{2\sigma_d R^d} \int_{\mathbb{B}_R} |x|^2\,\rd x +M\nonumber\\
	&= \frac{\kappa}{\lambda^2}\left(\int_{\mathbb{B}_{1}}\big(E_d*\zeta^{1/(m-1)}\big)(y)\,\rd y + \int_{\mathbb{B}_{\lambda R}\setminus \mathbb{B}_{1}}\big(E_d*\zeta^{1/(m-1)}\big)(y)\,\rd y\right) \nonumber\\
	& \qquad +\frac{R^2 M}{2(d+2)}+M \nonumber\\
	& \le \kappa \int_{\mathbb{B}_{1}}\big(E_d*\zeta^{1/(m-1)}\big)(y)\,\rd y + \frac{\kappa M_c}{(d-2) \sigma_d \lambda^2} \int_{\mathbb{B}_{\lambda R}\setminus \mathbb{B}_{1}} \frac{\rd y}{|y|^{d-2}} + (1+R^2)M \nonumber\\
	& \le  \kappa \int_{\mathbb{B}_{1}}\big(E_d*\zeta^{1/(m-1)}\big)(y)\,\rd y + \kappa R^2 M_c + (1+R^2)M =:\bar{M}_0\,. \label{q1}
\end{align}
Consequently,  recalling~\eqref{ub01}, we conclude that, for $\bar M\ge \bar{M}_0$,
\begin{equation}
    (u_{\lambda},w_{\lambda}) \in \mathcal{I}_{M,\bar{M},rad}(\mathbb{B}_R)\,. \label{ub00}
\end{equation}
Similarly, we compute
\begin{align*}
	\|v_\lambda\|_2^2&\le\frac{3\kappa^2}{\lambda^{4-d}} \int_{\mathbb{B}_{\lambda R}}\big|\big(E_d*\zeta^{1/(m-1)}\big)(y)\big|^2\,\rd y + \frac{3M^2}{4\sigma_d^2 R^{2d}} \int_{\mathbb{B}_R} |x|^4\,\rd x  +\frac{3M^2}{|\mathbb{B}_R|}\\
	&=\frac{3\kappa^2}{\lambda^{4-d}} \left( \int_{\mathbb{B}_{1}}\big|\big(E_d*\zeta^{1/(m-1)}\big)(y)\big|^2\,\rd y+\int_{\mathbb{B}_{\lambda R}\setminus \mathbb{B}_{1}}\big|\big(E_d*\zeta^{1/(m-1)}\big)(y)\big|^2\,\rd y \right)\\
	& \qquad + \frac{3R^{4-d} M^2}{4\sigma_d(d+4)}+\frac{3M^2}{|\mathbb{B}_R|}\,. 
\end{align*}
Invoking again~\eqref{eq.B3} we see that
\begin{align*}
	\int_{\mathbb{B}_{\lambda R}\setminus \mathbb{B}_{1}}\big|\big(E_d*\zeta^{1/(m-1)}\big)(y)\big|^2\,\rd y & = \frac{M_c^2}{(d-2)^2 \sigma_d} \int_1^{\lambda R} r^{3-d}\,\rd r \\
	& \le \frac{M_c^2}{(d-2)^2 \sigma_d} \left\{ \begin{array}{ll}
		\lambda R & (d=3)\,, \\
		\log{(\lambda R)} & (d=4)\,, \\
		\frac{1}{d-4} & (d\ge 5)\,,
	\end{array}\right.
\end{align*}
from which we deduce that
\begin{equation*}
	\int_{\mathbb{B}_{\lambda R}\setminus \mathbb{B}_{1}}\big|\big(E_d*\zeta^{1/(m-1)}\big)(y)\big|^2\,\rd y \le \frac{M_c^2(1+R)}{\sigma_d} \lambda\,.
\end{equation*}
Hence, there exists $c_1=c_1(d,R,M)>0$ such that
\begin{equation}\label{q2}
	\|v_\lambda\|_2^2\le c_1\big(1+\lambda^{d-3}\big)\,.
\end{equation}
Now, estimating the first three terms of the right-hand side of
\begin{align*}
    \mathcal{L}(u_{\lambda},w_{\lambda}) & = \int_{\mathbb{B}_R} \phi(u_{\lambda})(x)\,\rd x - \int_{\mathbb{B}_R} (u_{\lambda} v_{\lambda})(x)\,\rd x\\
&\quad + \frac{1}{2} \int_{\mathbb{B}_R} |\nabla v_{\lambda}(x)|^2\;\rd x + \frac{1}{2} \int_{\mathbb{B}_R} |v_{\lambda}(x)|^2\;\rd x
\end{align*}
 as in the proof of \cite[Proposition~4.4]{LW26} (see \cite[Equations~(4.16),(4.17),(4.19)]{LW26} therein with $V_\lambda=v_\lambda-M/|\mathbb{B}_R|$ instead) and using~\eqref{q2} for the last term, we infer the existence of constants $(c_2,c_3)\in (0,\infty)^2$ depending on $d$, $R$ and $M$ such that
\begin{equation*}
    \mathcal{L}(u_{\lambda},w_{\lambda}) \le - c_2 \lambda^{d-2} +c_3 + c_1\lambda^{d-3}\,.
\end{equation*}
Since $d\ge 3$, it readily follows from~\eqref{ub00} and the above inequality that
\begin{equation*}
    \mu_{M,\bar{M}} \le \liminf_{\lambda\to\infty} \mathcal{L}(u_{\lambda},w_{\lambda}) = -\infty\,,
\end{equation*}
and the proof is complete.
\end{proof}

We are now in a position to provide the proof of part~\textbf{(I)} of Theorem~\ref{THM1} which follows the lines of~\cite[Section~5]{LW26}:

\begin{proof}[Proof of Theorem~\ref{THM1}~\textbf{(I)}]
Let $M\in (0,M_c)$, $\bar{M}\ge \bar{M}_0$, and $(u^0,w^0)\in \mathcal{I}_{M,\bar{M},rad}(\mathbb{B}_R)$. Then the corresponding global solution $(u,w)$  to~\eqref{E} is such that $(u(t),w(t))$ belongs to $\mathcal{I}_{M,\bar{M},rad}(\mathbb{B}_R)$ for all $t\ge 0$ according to Theorem~\ref{T1:Ex}. Moreover, since
\begin{equation}\label{l9}
	\mathcal{L}(u,w)(t)  \le \mathcal{L}\big(u^0,w^0\big)\,, \quad t\ge 0\,,
\end{equation}
by Proposition~\ref{PP2}, the lower bound~\eqref{lbL} implies that $v\in L_\infty((0,\infty),W_2^1(\Omega))$. This  together with~\eqref{l9} and the properties of $\phi$ entails that  $u\in L_\infty((0,\infty),L_m(\Omega))$, so that Theorem~\ref{T1:Ex} guarantees that 
\begin{equation*}
	\sup_{t\ge 0}\big(\|u(t)\|_\infty+\|w(t)\|_\infty\big)<\infty\,.
\end{equation*}
We refer to \cite[Section~5]{LW26} for details.
\end{proof}

\section{Proof of Theorem~\ref{THM1}~\textbf{(II)}}\label{sec.4}

We turn to the proof of part~\textbf{(II)} of Theorem~\ref{THM1} and, as in \cite{Ho2002, HW2001}, we show that the energy functional $\mathcal{L}$ is bounded from below on the set of radially symmetric stationary solutions to~\eqref{E}. 
In the following, given $M>0$ we write $(u,w)\in\mathcal{S}_{M,rad}$ if and only if 
\begin{subequations}\label{Erad}
\begin{equation}
	u = w \in W_{p,rad}^{2,+}(\mathbb{B}_R)\,,\qquad \|u\|_1=M\,,
\end{equation}
and
\begin{align}
	\mathrm{div}\big(u\nabla(\phi'(u)-v)\big)&=0\ \text{ in }\ \mathbb{B}_R\,,\label{E1rad}\\
	\nabla u\cdot \textbf{n} & = 0 \ \text{ on }\ \partial\mathbb{B}_R \,,\label{NBCrad} \\
	v & = \mathcal{K}[u]\ \text{ in }\ \mathbb{B}_R\,,\label{E2rad}
\end{align}
\end{subequations}
recalling that $p\in (2d,\infty)$ is fixed.

\medskip

In fact, similarly as in~\cite{LW26}, the analysis of $\mathcal{S}_{M,rad}$ is captured by radially symmetric, positive solutions to a single elliptic equation. In the following, for a given $a>0$, we denote by $\Sigma_a$ the set of functions 
\begin{equation*}
	z\in W_{p,rad}^{2}(\mathbb{B}_R)\,, \quad z(x)>0 \;\text{ for }\; x\in \mathbb{B}_R\,,\quad  \langle z\rangle = a\,,  
\end{equation*}
solving
\begin{equation*}
	-\Delta \phi'(z)+\phi'(z)-\langle \phi'(z)\rangle=z-a\ \text{ in }\ \mathbb{B}_R\,,\qquad \nabla \phi'(z) \cdot \mathbf{n} = 0 \ \text{ on }\ \partial\mathbb{B}_R\,. 
\end{equation*}
Then $u\in \Sigma_{M/|\mathbb{B}_R|}$ for any $(u,w)\in \mathcal{S}_{M,rad}$ as shown next: 

\begin{prop}\label{L19}
Let $M>0$ and $(u,w)\in\mathcal{S}_{M,rad}$. Then $u\in\Sigma_a$ for $a:=M/|\mathbb{B}_R|$ and
\begin{equation}\label{n1}
	u(x)\ge\min\left\{1,\exp\left( \big(\langle \phi'(u)\rangle - a\big)/m \right)\right\}>0\,,\quad x\in\bar{\mathbb{B}}_R\,.
\end{equation}
Moreover,  with $v = \mathcal{K}[u]$ it holds that 
\begin{equation}\label{n22}
	\|v\|_1=M\,, \quad 0\le v=\phi'(u)-\langle \phi'(u)\rangle+a \le \|u\|_\infty \ \text{ in }\ \mathbb{B}_R\,.
\end{equation}
Finally,
\begin{equation}\label{n3}
	\|\nabla\phi'(u)\|_2^2+	\|\phi'(u)\|_2^2=|\mathbb{B}_R|\langle \phi'(u)\rangle^2+\int_\Omega \phi'(u)\left(u-a\right)\,\rd x\,.
\end{equation}
\end{prop}

\begin{proof}
The proof is similar to that of \cite[Proposition~6.1]{LW26}.
By premise, $u$ is continuous on $\bar{\mathbb{B}}_R$. Set $P:=\{x\in\mathbb{B}_R\,:\, u(x)>0\}$ and note from \eqref{E1rad} and~\eqref{NBCrad} that 
\begin{equation}\label{pp}
\nabla(\phi'(u)-v)=0\ \text{ in }\ P\,.
\end{equation}
Given any connected component $P_c$ of $P$, we thus find a constant $k_c\in\mathbb{R}$ such that 
\begin{equation}\label{n20}
	\phi'(u)=v+k_c\ \text{ in }\ P_c\,.
\end{equation}
Now, either $\min_{\bar{P}_c} u \ge 1$ so that
\begin{equation*}
	P_c=P=\mathbb{B}_R \;\;\text{ and }\;\; u(x) \ge 1\,, \quad x\in\bar{\mathbb{B}}_R\,.
\end{equation*}
Or $\min_{\bar{P}_c} u < 1$ and it follows from~\eqref{phi},~\eqref{n20}, and the non-negativity of $v=\mathcal{K}[v]$ that, for any $x\in P_c$ with $0< u(x)<1$,
\begin{equation*}
	-k_c\ge	-v(x)-k_c=-\phi'(u(x))=m\int_{u(x)}^1 \frac{(1+z)^{m-1}}{z}\,\rd z\ge -m\log(u(x))
\end{equation*}
and thus
\begin{equation*}
	u(x)\ge \exp\left( k_c/m \right)>0\,.
\end{equation*}
Thus, we conclude also in this case that
\begin{equation*}
	P_c=P=\mathbb{B}_R \;\;\text{ and }\;\; u(x) \ge \exp\left( k_c/m\right) >0\,, \quad x\in\bar{\mathbb{B}}_R\,.
\end{equation*}
Consequently, $P=\mathbb{B}_R$, and 
\begin{equation*}
	u(x)\ge\min\left\{1,\exp\left( k_c/m \right)\right\}>0\,,\quad x\in\bar{\mathbb{B}}_R\,.
\end{equation*}
Next, on the one hand, invoking~\eqref{E2rad} and~\eqref{n20} we deduce that
\begin{equation}\label{i1}
	-\Delta \phi'(u)=-\Delta v=u-v \ \text{ in }\ \mathbb{B}_R\,.
\end{equation}
In particular, the maximum principle implies that $v\le \|u\|_\infty$ in $\bar{\mathbb{B}}_R$, while the identity $\|v\|_1=\|u\|_1=M$ is an immediate consequence of~\eqref{E2rad} and~\eqref{n22}. On the other hand, integrating~\eqref{n20} over $\mathbb{B}_R$ and using~\eqref{int} yield 
\begin{equation}\label{i2}
	k_c=\langle\phi'(u)\rangle-\langle v\rangle =\langle\phi'(u)\rangle- \frac{M}{|\mathbb{B}_R|}\,.
\end{equation}
 Gathering~\eqref{n20},~\eqref{i1}, and ~\eqref{i2} we deduce that
\begin{equation}\label{i3}
	-\Delta \phi'(u)+\phi'(u)-\langle \phi'(u)\rangle=u-a\ \text{ in }\ \mathbb{B}_R\,,\qquad \nabla \phi'(u) \cdot \mathbf{n} = 0 \ \text{ on }\ \partial\mathbb{B}_R\,,
\end{equation}
with $a=M/|\mathbb{B}_R|$. Therefore,  $u\in \Sigma_a$ and we have shown~\eqref{n1} and~\eqref{n22}. Finally,~\eqref{n3} follows directly from~\eqref{i3}.
\end{proof}

\begin{cor}\label{C20}
Let $M>0$ and consider $(u,w)\in\mathcal{S}_{M,rad}$. Then
\begin{equation}\label{n4}
	mathcal{L}(u,w) = \mathcal{L}(u,u) =\int_{\mathbb{B}_R} \phi(u)\,\rd x -\frac{1}{2} \int_{\mathbb{B}_R} u\phi'(u)\,\rd x+\frac{M}{2}\langle \phi'(u)\rangle-\frac{M^2}{2|\mathbb{B}_R|}\,.
\end{equation}
Moreover, there exists a constant $\ell(M)>0$ such that
\begin{equation}\label{n5}
	\mathcal{L}(u,u)\ge \ell(M)+\frac{M}{2}\langle \phi'(u)\rangle\,.
\end{equation}
\end{cor}

\begin{proof}
With $v=\mathcal{K}[u]$ we infer from~\eqref{E2rad} that
\begin{equation*}
	\int_{\mathbb{B}_R} uv\,\rd x = \int_{\mathbb{B}_R} \big(v-\Delta v\big) v\,\rd x =  \int_{\mathbb{B}_R} |v|^2\,\rd x+ \int_{\mathbb{B}_R} |\nabla v|^2\,\rd x\,.
\end{equation*}
Therefore, using also~\eqref{n22} (with $a=M/|\mathbb{B}_R|$), it follows that
\begin{align*}
	\mathcal{L}(u,u) & = \int_{\mathbb{B}_R} \phi(u)\,\rd x - \int_{\mathbb{B}_R} uv\,\rd x + \frac{1}{2} \int_{\mathbb{B}_R} |\nabla v|^2\,\rd x+ \frac{1}{2}\int_{\mathbb{B}_R} | v|^2\,\rd x\\
	& = \int_{\mathbb{B}_R} \phi(u)\,\rd x - \frac{1}{2} \int_{\mathbb{B}_R} |\nabla v|^2\,\rd x- \frac{1}{2}\int_{\mathbb{B}_R} | v|^2\,\rd x\\
&=\int_{\mathbb{B}_R} \phi(u)\,\rd x -\frac{1}{2} \int_{\mathbb{B}_R} |\nabla \phi'(u)|^2\,\rd x-\frac{1}{2}\int_{\mathbb{B}_R} |\phi'(u)-\langle \phi'(u)\rangle+a|^2\,\rd x\,.
\end{align*}
We use~\eqref{i3} to derive the identity
\begin{align*}
\int_{\mathbb{B}_R} |\phi'(u)-\langle \phi'(u)\rangle+a|^2\,\rd x&=\int_{\mathbb{B}_R} \big(\phi'(u)-\langle \phi'(u)\rangle+a\big)\big(u+\Delta \phi'(u))\big)\,\rd x\\
&= \int_{\mathbb{B}_R} u\phi'(u)\,\rd x -M\langle\phi'(u)\rangle +Ma- \int_{\mathbb{B}_R} |\nabla\phi'(u)|^2\,\rd x \,,
\end{align*}
and we combine the above two formulas to obtain~\eqref{n4}. Finally, since Lemma~\ref{LTJ} bounds the function $r \mapsto 2\phi(r)-r\phi'(r)$ on $(0,\infty)$ from below,~\eqref{n5} is implied by~\eqref{n4}.
\end{proof}

As in \cite[Lemma~6.3]{LW26}, we derive additional estimates for $u\in \Sigma_a$. We emphasize that different arguments are needed at several places due to the zero-order linear term in the equation determining $v$.

\begin{lem}\label{Lx}
Let $a>0$, $M=a|\mathbb{B}_R|$ and consider $u\in \Sigma_a$. Then, recalling that 
\begin{equation*}
    v=\mathcal{K}[u] = \phi'(u)-\langle \phi'(u)\rangle +a\,,
\end{equation*}
the functions $u$, $v$, and $\phi'(u)$ satisfy 
\begin{equation}\label{n8}
	|\partial_r v(r)| = |\partial_r \phi'(u(r))| \le \min\left\{\frac{ \|u\|_\infty }{d} r, \frac{M}{\sigma_d}r^{1-d}\right\}\,,\quad r\in (0,R]\,,
\end{equation}
and
\begin{equation}\label{n9}
	|\phi'(u(r))-\langle\phi'(u)\rangle| \le \frac{M}{\sigma_d(m-1)}\left(\frac{r^{2-d}}{d}+\frac{R^{2-d}}{2}\right)\,,\quad r\in (0,R]\,.
\end{equation}
In particular,
\begin{equation}\label{n7}
	|\partial_r u(r)| \le \min\{u(r),u(r)^{2-m}\}\|u\|_\infty\frac{r}{dm}\,,\quad r\in (0,R]\,,
\end{equation}
and 
\begin{equation}\label{n7x}
	|\partial_r u(r)| \le \frac{M}{m\sigma_d}\min\{u(r),u(r)^{2-m}\}r^{1-d}\,,\quad r\in (0,R]\,.
\end{equation}
Moreover, there is $c_1=c_1(d,R,M)>0$ such that
\begin{align}\label{N7A}
	\partial_r u(r)\le c_1 \min\{u(r),u(r)^{2-m}\} r^{3-d}\,,\quad r\in (0,R]\,,
\end{align}
and 
\begin{equation}\label{N6}
	0<u(r)\le \left[ \frac{dM}{\sigma_d} + c_1 \left( \frac{M}{\sigma_d} \right)^{2/d} R^2 \right] r^{-d}\,,\quad r\in (0,R]\,.
\end{equation}
\end{lem}

\begin{proof} 
Since $v=\mathcal{K}[u]$ is radially symmetric, we have
\begin{equation*}
	-\partial_r\left(r^{d-1}\partial_r v(r)\right)= r^{d-1} (u-v)(r)\,,\quad r\in (0,R]\,,
\end{equation*}
so that integration yields
\begin{align}\label{g1}
	-r^{d-1}\partial_r v(r) = \int_0^r s^{d-1} (u-v)(s)\,\rd s \,,\quad r\in [0,R] \,.
\end{align}
On the one hand, since $v$ is positive, we obtain from~\eqref{eqi} and~\eqref{g1} that
\begin{align*}
	-r^{d-1}\partial_r v(r) \le \int_0^r s^{d-1} u(s)\,\rd s\le  \min\left\{\frac{ \|u\|_\infty }{d} r^d, \frac{M}{\sigma_d}\right\} \,,\quad r\in [0,R] \,;
\end{align*}
that is,
\begin{align}\label{i99}
- \partial_r v(r) \le\min\left\{\frac{ \|u\|_\infty }{d} r, \frac{M}{\sigma_d}r^{1-d}\right\}\,,\quad r\in [0,R] \,.
\end{align}
On the other hand, since $u$ is positive, we derive from~\eqref{eqi} and~\eqref{g1} that
\begin{equation*}
	-r^{d-1}\partial_r v(r) \ge -\int_0^r s^{d-1} v(s)\,\rd s \ge - \min\left\{ \frac{\|v\|_1}{\sigma_d} , \frac{\|v\|_\infty}{d} r^d \right\}\,, \quad r\in (0,R]\,,
\end{equation*}
hence, using~\eqref{n22},
\begin{equation}\label{po}
	\partial_r v(r)\le \min\left\{ \frac{M}{\sigma_d}  r^{1-d}, \frac{\|u\|_\infty}{d} r\right\}\,,\quad r\in (0,R]\,.
\end{equation} 
Gathering~\eqref{i99} and~\eqref{po}, we deduce~\eqref{n8} in view of $\partial_r v=\partial_r \phi'(u)$ by~\eqref{n22}. Furthermore, since
\begin{equation}\label{g2}
	\partial_r\phi'(u(r))=\phi''(u(r))\partial_r u(r)=m\frac{(1+u(r))^{m-1}}{u(r)}\partial_r u(r)\,,\quad r\in (0,R]\,,
\end{equation}
we derive~\eqref{n7} and~\eqref{n7x} from~\eqref{n8} and the positivity of $u$. 

Next, as in \cite[Lemma~6.3]{LW26}, we use~\eqref{n8} to obtain, for $r\in (0,R]$,
\begin{align}
	|\phi'(u(r))-\langle\phi'(u)\rangle| &= \left| \frac{d}{R^d}\int_0^R \big(\phi'(u)(r)-\phi'(u)(s)\big)\, s^{d-1}\,\rd s \right| \nonumber\\
	&= \left| \frac{d}{R^d}\int_0^R\int_s^r \partial_r\phi'(u)(s_*)\,\rd s_*\, s^{d-1}\,\rd s \right| \label{ccc}\\
	&\le \frac{dM}{R^d\sigma_d} \int_0^R \left| \int_s^r s_*^{1-d}\,\rd s_*\right| \, s^{d-1}\,\rd s \nonumber\\
	&\le  \frac{dM}{R^d\sigma_d(d-2)} \int_0^R\left(  r^{2-d}+s^{2-d} \right)\, s^{d-1}\,\rd s \nonumber\\
	&= \frac{dM}{R^d\sigma_d(d-2)} \left(\frac{r^{2-d}R^d}{d}+\frac{R^2}{2}\right)\,, \nonumber
\end{align}
which is~\eqref{n9}. Now, for $r\in (0,R]$, we infer from~\eqref{n22},~\eqref{n9},~\eqref{g1}, and the positivity of~$u$ that
\begin{align*} 
	-r^{d-1}\partial_r v(r) &\ge -\int_0^r s^{d-1} |\phi'(u)(s)-\langle \phi'(u)\rangle|\,\rd s -\frac{ar^d}{d}\\
    &\ge -\frac{M}{\sigma_d(m-1)}\int_0^r s^{d-1} \left(\frac{s^{2-d}}{d}+\frac{R^{2-d}}{2}\right)\,\rd s-\frac{ar^d}{d} \\
    & = - \frac{M}{\sigma_d (m-1)} \left[ \frac{r^2}{2d} + \frac{R^{2-d} r^d}{2d} \right] - \frac{a r^d}{d} 
\end{align*}
and hence
\begin{equation}\label{615} 
	\partial_r\phi'(u)(r) =\partial_r v(r) \le m c_1 r^{3-d}\,,\quad r\in (0,R]\,,
\end{equation}
with
\begin{equation*}
    c_1 :=  \frac{1}{md} \left[ \frac{M}{\sigma_d}(m-1) + a R^{d-2} \right]\,.
\end{equation*}
Combining~\eqref{g2} and~\eqref{615} gives~\eqref{N7A}. 

Finally,~\eqref{N7A}, H\"older's inequality, and the positivity of $u$ imply that
\begin{align*}
	\int_0^r s^{d}\partial_r u(s)\,\rd s & \le c_1 \int_0^r s^3 u(s)^{2-m} \,\rd s = c_1 \int_0^r s^{(d+2)/d} \big(s^{d-1} u(s)\big)^{2/d}\,\rd s\\
    &\le c_1 \left( \int_0^r s^{d-1} u(s)\,\rd s \right)^{2/d} \left( \int_0^r s^{(d+2)/(d-2)}\,\rd s \right)^{(d-2)/d} \\ 
    & \le c_1 \left( \frac{M}{\sigma_d} \right)^{2/d} \left( \frac{(d-2)}{2d} r^{2d/(d-2)} \right)^{(d-2)/d} \\
    & \le c_1 \left( \frac{M}{\sigma_d} \right)^{2/d} r^2 \le c_1 \left( \frac{M}{\sigma_d} \right)^{2/d} R^2\,.
\end{align*}
Hence, writing
\begin{equation*}
	r^{d} u(r) =\int_0^r\left(d s^{d-1}u(s)+s^{d}\partial_r u(s)\right)\, \rd s \le \frac{dM}{\sigma_d} + c_1 \left( \frac{M}{\sigma_d} \right)^{2/d} R^2
\end{equation*}
for $r\in (0,R]$, we obtain~\eqref{N6}. This proves the lemma.
\end{proof}

 Recalling that $\zeta$ is given in~\eqref{b0}, we now prove the following result, which is fundamental to our purposes.

\begin{prop}\label{P20}
Let $a>0$, $M=a|\mathbb{B}_R|$, and consider a sequence $(u_l)_{l\ge 1}$ in $\Sigma_a$  such that $\lambda_l:=\|u_l\|_\infty^{1/d}\to \infty$ as $l\to\infty$. For $l\ge 1$, set 
\begin{equation}\label{t77}
	P_l(x):=\frac{1}{\lambda_l^d}u_l\left(\frac{x}{\lambda_l}\right)\,,\quad x\in \bar{\mathbb{B}}_{\lambda_l R}\,,
\end{equation} 
and
\begin{equation}\label{T7}
	\varphi_l(x):=\frac{1}{\lambda_l^{d-2}}\phi'\left(u_l\right)\left(\frac{x}{\lambda_l}\right)\,,\quad x\in \bar{\mathbb{B}}_{\lambda_l R}\,.
\end{equation}
Then there are $P\in C(\mathbb{R}^d)$ and $\varphi\in C^2(\mathbb{R}^d)$ such that, up to a subsequence, 
\begin{equation*}
	P_l\rightarrow P \ \text{ in }\ C(\mathbb{R}^d)\,,\qquad \varphi_l\rightarrow \varphi \ \text{ in }\ C^2(\mathbb{R}^d)\,.
\end{equation*}
Both $P$  and $\varphi$ are non-increasing with respect to $r=|x|$ with $P(0)=1$, and $P$ has compact support $\bar{\mathbb{B}}_{\rho_0}$ for some $\rho_0\in (0,\infty)$. Furthermore,  $\varphi$ satisfies
\begin{equation}\label{varphi3}
	-\Delta \varphi=\left(\frac{m-1}{m}\right)^{1/(m-1)}\varphi_+^{1/(m-1)}\ \text{ in }\ \mathbb{R}^d
\end{equation}
with $\varphi_+=\max\{\varphi,0\}$, 
\begin{subequations}\label{420}
\begin{equation}
	P^{m-1} = \dfrac{m-1}{m} \varphi_+ \;\;\text{ in }\;\; \mathbb{R}^d\,, \label{420a}
\end{equation} 
and
\begin{equation}\label{420b}
 	\varphi(x)=\left\{\begin{array}{ll}
	\dfrac{m}{m-1} \rho_0^{2-d} \zeta\left(\dfrac{x}{\rho_0}\right)\,, & x\in \bar{\mathbb{B}}_{\rho_0}\,,\\[20pt]
	\dfrac{M_c}{\sigma_d}\dfrac{|x|^{2-d}-\rho_0^{2-d}}{d-2} \,, & x\in \R^d\setminus\bar{\mathbb{B}}_{\rho_0}\,.
\end{array}\right.
\end{equation}
\end{subequations}
Finally, it holds that $M\ge M_c$ with critical mass $M_c$ defined in~\eqref{Mc}.
\end{prop}

\begin{proof}
We divide the proof into several steps. We begin with estimates on $(P_l)_{l\ge 1}$ and go on with estimates on $(\varphi_l)_{l\ge 1}$. The main outcome of the first two steps is then the compactness of $(P_l)_{l\ge 1}$ and $(\varphi_l)_{l\ge 1}$ in suitable spaces and thus the existence of cluster points $(P,\varphi)$ as $l\to\infty$. The last step is to identify $P$ and $\varphi$, eventually leading us to~\eqref{420}. 

\smallskip

\noindent\textbf{Step~1: Pointwise estimates on $(P_l)_{l\ge 1}$}.
Let $l\ge 1$. Since $u_l$ belongs to $\Sigma_a$, we infer from~\eqref{n1} and~\eqref{t77} that 
\begin{equation}\label{n46}
	0< P_l(x)\le \|P_l\|_{L_\infty(\mathbb{B}_{\lambda_l R})} = 1\,,\quad x\in \bar{\mathbb{B}}_{\lambda_l R}\,,\quad \text{and}\quad \|P_l\|_{L_1(\mathbb{B}_{\lambda_l R})} = M\,,
\end{equation} 
and note from~\eqref{n7} that
\begin{equation}\label{n50}
 	|\partial_r P_l(r)| = \frac{1}{\lambda_l^{d+1}} \left| \partial_r u_l\left(\frac{r}{\lambda_l}\right)\right| \le \frac{1}{\lambda_l^{d+2}} \left[ u_l\left(\frac{r}{\lambda_l}\right) \right]^{2-m} \|u_l\|_\infty \frac{r}{dm}\le\frac{r}{dm} 
\end{equation}
for $r\in [0,\lambda_l R]$, since $\lambda_l=\|u_l\|_\infty^{1/d}$ and $m=2-2/d$. Next, by~\eqref{phi} and~\eqref{t77},
\begin{align}
	\partial_r\big(\lambda_l^{-d}+P_l(r)\big)^{m} &= \frac{m}{\lambda_l^{d+1}} \left(\lambda_l^{-d}+\lambda_l^{-d}u_l\left(\frac{r}{\lambda_l}\right)\right)^{m-1} \partial_r u_l\left(\frac{r}{\lambda_l}\right)\nonumber\\
	&=\frac{1}{\lambda_l^{2d-1}} \phi''\left(u_l\right)\left(\frac{r}{\lambda_l}\right) u_l\left(\frac{r}{\lambda_l}\right) \partial_ru_l\left(\frac{r}{\lambda_l}\right)\nonumber\\
	&= \frac{1}{\lambda_l^{d-1}} \partial_r\phi'\left(u_l\right) \left(\frac{r}{\lambda_l}\right) P_l(r)\,. \label{n50b}
\end{align}
Thus, we may use~\eqref{615} and~\eqref{n46} to obtain
\begin{equation}
	\partial_r\big(\lambda_l^{-d}+P_l(r)\big)^{m} \le \frac{1}{\lambda_l^{d-1}} \frac{m c_1 r^{3-d}}{\lambda_l^{3-d}} = \frac{m c_1 r^{3-d}}{\lambda_l^{2}}\,. \label{36}
\end{equation}
Integrating~\eqref{36} gives, for $0< r_1\le r_2\le \lambda_l R$,
\begin{equation}\label{Ps1}
 	\big(\lambda_l^{-d}+P_l(r_2)\big)^{m}- \big(\lambda_l^{-d}+P_l(r_1)\big)^{m}\le \frac{m c_1}{\lambda_l^2} \int_{r_1}^{r_2}s^{3-d}\,\rd s\,.
\end{equation}
Fix $\epsilon\in  (0,1)$ and $0<\theta<(2-\epsilon)/(d+\epsilon-4)_+$. On the one hand, if $\lambda_l^{-\theta}<r_2\le\lambda_l R$, then~\eqref{Ps1} gives 
\begin{align*}
	\big(\lambda_l^{-d} + P_l(r_2)\big)^{m} - \big(\lambda_l^{-d}+P_l(\lambda_l^{-\theta})\big)^{m} & \le \frac{m c_1 r_2^\epsilon}{\lambda_l^2} \int_{\lambda_l^{-\theta}}^{r_2}s^{3-\epsilon -d}\,\rd s \\
	& \le \frac{m c_1 R^\epsilon}{\lambda_l^{2-\epsilon}} \frac{\lambda_l^{-\theta(4-\epsilon -d)}-r_2^{4-\epsilon -d}}{d+\epsilon-4}\,.
\end{align*} 
Now, either $d=3$ and 
\begin{align*}
	\frac{m c_1 R^\epsilon}{\lambda_l^{2-\epsilon}}\frac{\lambda_l^{-\theta(4-\epsilon -d)}-r_2^{4-\epsilon -d}}{d+\epsilon-4} & = \frac{m c_1 R^\epsilon}{\lambda_l^{2-\epsilon}} \frac{r_2^{1-\epsilon}-\lambda_l^{-\theta(1-\epsilon )}}{1-\epsilon} \\
	& \le \frac{m c_1 R^\epsilon}{\lambda_l^{2-\epsilon}}\frac{(\lambda_l R)^{1-\epsilon}}{1-\epsilon} = \frac{m c_1 R}{(1-\epsilon)\lambda_l}\,.
\end{align*}
Or $d\ge 4$ and
\begin{equation*}
	\frac{m c_1 R^\epsilon}{\lambda_l^{2-\epsilon}}\frac{\lambda_l^{-\theta(4-\epsilon -d)}-r_2^{4-\epsilon -d}}{d+\epsilon-4}\le \frac{m c_1 R^\epsilon}{\lambda_l^{2-\epsilon}}\frac{\lambda_l^{-\theta(4-\epsilon -d)}}{d+\epsilon-4} = \frac{m c_1 R^\epsilon}{(d+\epsilon-4)\lambda_l^{2-\epsilon+\theta(4-\epsilon -d)}}\,.
\end{equation*}
Setting $\eta:=1$ if $d=3$ and $\eta:=2-\epsilon+\theta(4-\epsilon -d)>0$ if $d\ge 4$, we conclude that
\begin{align}\label{opp}
	\big(\lambda_l^{-d} + P_l(r_2)\big)^{m} \le \big(\lambda_l^{-d}+P_l(\lambda_l^{-\theta})\big)^{m} +  \frac{m c_1(R+R^\epsilon)}{|d+\epsilon-4|\lambda_l^\eta}\,,\quad \lambda_l^{-\theta}<r_2\le\lambda_l R\,.
\end{align} 
On the other hand, if $0<r_2\le \lambda_l^{-\theta}$, it follows from~\eqref{n46} and~\eqref{n50} that
\begin{align*}
	\big| \big(\lambda_l^{-d} + P_l(r_2)\big)^{m} & -\big(\lambda_l^{-d} + P_l(\lambda_l^{-\theta})\big)^{m} \big| \le m\int_{r_2}^{\lambda_l^{-\theta}} \big(\lambda_l^{-d} + P_l(s)\big)^{m-1} \left|\partial_r P_l(s)\right|\,\rd s\\
	&\quad \le \frac{1}{d} \big(\lambda_l^{-d} +1\big)^{m-1} \int_{r_2}^{\lambda_l^{-\theta}} s\,\rd s \le  \frac{\big(\lambda_l^{-d} +1\big)^{m-1}}{2d \lambda_l^{2\theta}}\,.
\end{align*}
Combining this inequality with~\eqref{opp} yields
\begin{align*}
    \big(\lambda_l^{-d}+P_l(r_2)\big)^{m} &\le \big(\lambda_l^{-d} + P_l(\lambda_l^{-\theta})\big)^{m} + \frac{\big(\lambda_l^{-d} +1\big)^{m-1}}{2d \lambda_l^{2\theta}} + \frac{m c_1 (R+R^\epsilon)}{|d+\epsilon-4|\lambda_l^\eta}
\end{align*}
for $0<r_2\le \lambda_l R$. Consequently,
\begin{align*}
     \big(\lambda_l^{-d}+1\big)^{m}&= \big(\lambda_l^{-d}+\|P_l\|_{L_\infty(\mathbb{B}_{\lambda_l R})}\big)^{m}\\
	& \le \big(\lambda_l^{-d} + P_l(\lambda_l^{-\theta})\big)^{m} + \frac{\big(\lambda_l^{-d} +1\big)^{m-1}}{2d \lambda_l^{2\theta}} + \frac{m c_1(R+R^\epsilon)}{|d+\epsilon-4|\lambda_l^\eta}\,,
\end{align*}
so that the subadditivity of $z\mapsto z^{1/m}$ implies that
\begin{equation*}
    1\le P_l(\lambda_l^{-\theta})+\left( \frac{\big(\lambda_l^{-d} +1\big)^{m-1}}{2d \lambda_l^{2\theta}} + \frac{m c_1(R+R^\epsilon)}{|d+\epsilon-4|\lambda_l^\eta} \right)^{1/m}\,.
\end{equation*}
Finally, due to~\eqref{n50} we have
\begin{align*}
	P_l(\lambda_l^{-\theta}) & \le P_l(0)+\int_0^{\lambda_l^{-\theta}}\partial_r P_l(s)\,\rd s \\
	& \le P_l(0)+\frac{1}{md}\int_0^{\lambda_l^{-\theta}}s\,\rd s = P_l(0) + \frac{1}{2md \lambda_l^{2\theta}}\,.
\end{align*}
Combining the previous two inequalities leads us to
\begin{equation}\label{6.x}
    1\le P_l(0)+Q(\lambda_l)\,,
\end{equation}
where
\begin{equation*}
	Q(\lambda_l):= \frac{1}{2md \lambda_l^{2\theta}} + \left( \frac{\big(\lambda_l^{-d} +1\big)^{m-1}}{2d\lambda_l^{2\theta}} + \frac{m c_1(R+R^\epsilon)}{|d+\epsilon-4|\lambda_l^\eta}\right)^{1/m}\,.
\end{equation*}
Let us observe at this point that 
\begin{equation}\label{popp}
\lim_{l\to\infty} Q(\lambda_l)= 0\,,
\end{equation}
since $\theta>0$, $\eta>0$, and $d\ge 3$.

\medskip
\noindent\textbf{Step~2: Pointwise estimates on $(\varphi_l)_{l\ge 1}$.} 
Let $l\ge 1$. We  proceed as in the proof of~\cite[Proposition~6.4]{LW26} and use~\eqref{phi},~\eqref{t77}, and~\eqref{T7} to derive
\begin{align}
	\varphi_l(r) & = \frac{1}{\lambda_l^{d-2}} \int_{1}^{u_l(r/\lambda_l)} \phi''(s)\,\rd s = \frac{m}{\lambda_l^{d-2}} \int_{\lambda_l^{-d}}^{P_l(r)} \frac{\left(1+\lambda_l^d s \right)^{m-1}}{s}\,\rd s \nonumber \\
	& =m \int_{\lambda_l^{-d}}^{P_l(r)} \frac{\left(\lambda_l^{-d}+s\right)^{m-1}}{s}\,\rd s\,,\quad r\in [0,\lambda_l R]\,. \label{6.y}
\end{align}
On the one hand, it follows from~\eqref{6.x} and~\eqref{popp} that there exists $l_1\in\N$ such that
\begin{equation*}
	\lambda_l^{-d}\le 1- Q(\lambda_l)\le P_l(0)\,,\quad l\ge l_1\,,
\end{equation*}
so that~\eqref{6.y} yields
\begin{align}
	\varphi_l(0) & \ge m \int_{\lambda_l^{-d}}^{P_l(0)}s^{m-2} \,\rd s=\frac{m}{m-1}\big(P_l(0)^{m-1}-\lambda_l^{-d(m-1)}\big)\nonumber\\
	& \ge \frac{m}{m-1} \left(\left[1- Q(\lambda_l)\right]^{m-1}-\lambda_l^{-d(m-1)}\right)\,,\label{6.zl}
\end{align}
for $l\ge l_1$. On the other hand, since $\lambda_l^{-d} \le 1$ for $l\ge l_1$, we infer from~\eqref{6.y} and the subadditivity of $z\mapsto z^{m-1}$ that, for $r\in [0,\lambda_l R]$ and $l\ge l_1$,
\begin{align*}
	\varphi_l(r)& \le m \int_{\lambda_l^{-d}}^{1} \left(\frac{\lambda_l^{-d(m-1)}}{s}+s^{m-2}\right) \,\rd s \\
	& =\frac{m}{\lambda_l^{d-2}}\left(-\log\lambda_l^{-d}\right)+\frac{m}{m-1}\left(1-\lambda_l^{-d(m-1)}\right)
\end{align*}
and thus
\begin{equation}\label{6.zu}
	\varphi_l(r) \le \frac{m}{m-1}\left(1-\frac{1}{{\lambda_l^{d-2}}}\right)+\frac{md\log\lambda_l}{\lambda_l^{d-2}}\,,\quad r\in [0,\lambda_l R]\,, \quad l\ge l_1\,.
\end{equation}
Making $l_1>0$ larger, if necessary, we deduce from~\eqref{popp},~\eqref{6.zl}, and~\eqref{6.zu} (with $r=0$) that
\begin{equation}\label{6.zz}
	\frac{m}{2(m-1)}\le \varphi_l(0)\le \frac{2m}{m-1}\,,\quad l\ge l_1 \,.
\end{equation}
Furthermore, it follows from~\eqref{n8} that, for $r\in [0,\lambda_l R]$,
\begin{align}
	|\partial_r\varphi_l(r)| & = \frac{1}{\lambda_l^{d-1}} \left| \partial_r\phi'(u_l)\left( \frac{r}{\lambda_l} \right) \right| \nonumber \\
	& \le \frac{1}{\lambda_l^{d-1}} \min\left\{\frac{ \|u_l\|_\infty }{d}  \frac{r}{\lambda_l}, \frac{M}{\sigma_d}\left( \frac{r}{\lambda_l}\right)^{1-d}\right\}=\min\left\{\frac{r}{d}  , \frac{Mr^{1-d}}{\sigma_d}\right\}\,. \label{x03}
\end{align}
Therefore,
\begin{align*}
|\varphi_l(r)|&\le\left|\varphi_l(0)\right|+\int_0^r\left|\partial_r\varphi_l(s)\right|\,\rd s\\
	& \le \left|\varphi_l(0)\right|+\int_0^{\left(Md/\sigma_d\right)^{1/d}}\frac{s}{d}\,\rd s +\int_{\left(Md/\sigma_d\right)^{1/d}}^\infty\frac{Ms^{1-d}}{\sigma_d}\,\rd s\,,
\end{align*}
which, combined   with~\eqref{6.zz}, gives 
\begin{equation}\label{6.34}
	|\varphi_l(r)|\le  \frac{2m}{m-1}+\frac{1}{2d}\left(\frac{Md}{\sigma_d}\right)^{2/d}+\frac{M}{\sigma_d (d-2)}\left(\frac{Md}{\sigma_d}\right)^{(2-d)/d} =:c_2
\end{equation}
for $r\in [0,\lambda_l R]$ and $l\ge l_1 $. Also note that
\begin{align*}
	\frac{1}{\lambda_l^{d-2}}\langle\phi'(u_l)\rangle &=\frac{d}{R^d\lambda_l^{d-2}}\int_0^Rr^{d-1}\phi'(u_l)(r)\,\rd r =  \frac{d}{(\lambda_l R)}^{d} \int_0^{\lambda_l R}r^{d-1}\varphi_l(r)\,\rd r\,,
\end{align*}
so that~\eqref{6.34} entails that
\begin{align}\label{54}
	\frac{1}{\lambda_l^{d-2}}\left|\langle\phi'(u_l)\rangle\right|&\le  c_2\,,\quad l\ge l_1\,.
\end{align}
By~\eqref{i3}, 
\begin{equation}
	-\Delta\varphi_l = S_l := - \frac{\varphi_l}{\lambda_l^2} + \frac{\langle\phi'(u_l)\rangle}{\lambda^d} + P_l - \frac{a}{\lambda_l^d} \ \text{ in }\  \mathbb{B}_{\lambda_l R}\,, \label{x01}
\end{equation}
and we infer from~\eqref{n46},~\eqref{6.34} and~\eqref{54} that, for $l\ge l_1$,
\begin{equation}
	\|S_l\|_{L_\infty(\mathbb{B}_{\lambda_l R})} \le \frac{2c_2}{\lambda_l^2} + 1 + \frac{a}{\lambda_l^d} \le 1+a+2c_2. \label{x02}
\end{equation}
Now, let $R_0>0$ be arbitrary. Then, there exists $l_2(R_0)\ge l_1$ such that $\lambda_l R\ge 3 R_0$ for $l\ge l_2(R_0)$. Thus, since~\eqref{6.34} ensures that 
\begin{equation*}
	\|\varphi_l\|_{L_\infty(\mathbb{B}_{3R_0})} \le c_2\,,\quad l\ge l_2(R_0)\,,
\end{equation*}
it follows from~\eqref{x01},~\eqref{x02}, and \cite[Theorem~9.11]{GT2001} that $(\varphi_l)_{l\ge l_2(R_0)}$ is bounded in $W_q^2(\mathbb{B}_{2R_0})$ for any $q\in (1,\infty)$. In particular, since $W_q^2(\mathbb{B}_{2R_0})$ is continuously embedded in $W_\infty^1(\mathbb{B}_{2R_0})$ for $q>d$, we deduce from~\eqref{n46},~\eqref{n50} \eqref{x03},~\eqref{6.34} and~\eqref{x02} that $(S_l)_{l\ge l_2(R_0)}$ is bounded in $W_\infty^1(\mathbb{B}_{2R_0})$. Applying \cite[Theorems~9.19]{GT2001} to~\eqref{x01} entails that, for any $q\in (1,\infty)$, there exists $c(q,R_0)>0$ such that
\begin{equation}\label{622}
	\left\|\varphi_l\right\|_{W_q^3(\mathbb{B}_{R_0})}\le c(q,R_0)\,,\quad l\ge l_2(R_0)\,.
\end{equation}

\medskip
\noindent\textbf{Step~3: Convergence.} 
Let $q\in (d,\infty)$. Owing to~\eqref{n46},~\eqref{n50} and~\eqref{622}, the sequences $(P_l)_{l\ge l_2(R_0)}$ and $(\varphi_l)_{l\ge l_2(R_0)}$ are bounded in $W_\infty^1(\mathbb{B}_{R_0})$ and $C^{3-(d/q)}(\mathbb{B}_{R_0})$, respectively, for each $R_0>0$. The Arzel\`a-Ascoli theorem, along with a diagonal process, then ensures the existence of functions $P\in C(\mathbb{R}^d)$ and $\varphi\in C^2(\mathbb{R}^d)$ such that, up to a subsequence,
\begin{equation}
		P_l\to P \ \text { in }\ C(\mathbb{R}^d)\,, \quad \varphi_l\to \varphi \ \text { in }\ C^2(\mathbb{R}^d)\,, \label{N51}
	\end{equation}
as $l\to \infty$ (i.e. locally uniformly in $\mathbb{R}^d$). Moreover,~\eqref{54} guarantees that there exists $\alpha\in\R$ such that, again up to a subsequence, 
\begin{equation}\label{N51x}
	\lim_{l\to\infty} \frac{\langle\phi'(u_l)\rangle}{\lambda_l^{d-2}} =\alpha\,.
\end{equation}

\medskip

\noindent\textbf{Step~4: Identifying $P$ and $\varphi$.} 
To begin with, let us recall that~\eqref{n50b} also reads
\begin{equation*}
	\partial_r P_l(r) = \frac{1}{m} \frac{P_l(r)}{\left(\lambda_l^{-d}+P_l(r)\right)^{m-1}} \partial_r\varphi_l(r)\,.
\end{equation*}
We then readily deduce from~\eqref{N51} that $P$ belongs to $C^1(\R^d)$ and satisfies
\begin{equation}\label{45}
	\partial_r P(r) = \frac{1}{m}P(r)^{2-m}\partial_r\varphi(r)\,,\quad r\ge 0\,,
\end{equation}
Moreover, combining~\eqref{x01},~\eqref{N51}, and ~\eqref{N51x} implies that $\varphi$ is a classical solution to 
\begin{equation}\label{t4}
	-\Delta \varphi=P \ \text{ in }\ \R^d\,.
\end{equation}
Next,~\eqref{Ps1} entails that $P^m$, and hence also $P$, is non-increasing. The latter, along with~\eqref{45}, implies that also $\varphi$ is non-increasing.  Moreover,~\eqref{n46},~\eqref{6.x},~\eqref{popp} and~\eqref{N51} yield $P(0)=1$, while~\eqref{6.zl},~\eqref{6.zu} and~\eqref{N51} ensure that 
\begin{equation*}
\varphi(0)=m/(m-1)\,.
\end{equation*}

Now, either $P(r)>0$ for every $r\in [0,\infty)$, or there is $\rho_0\in (0,\infty)$ such that $\mathrm{supp}(P)=\bar{\mathbb{B}}_{\rho_0}$. In fact, for each $r>0$ such that $P(r)>0$, we infer from~\eqref{6.y} and~\eqref{N51} that
\begin{equation}\label{6.yd}
    \varphi(r)=m \int_{0}^{P(r)}s^{m-2}\,\rd s=\frac{m}{m-1}P^{m-1}(r)\,.
\end{equation}
Assume for contradiction that $P>0$ on~$\mathbb{R}^d$. Then~\eqref{t4} and~\eqref{6.yd} imply that
$\varphi$ is a positive (classical) solution to
\begin{equation*}
	-\Delta \varphi=\left(\frac{m-1}{m}\right)^{1/(m-1)} \varphi^{1/(m-1)}\ \text{ in }\ \R^d\,,
\end{equation*}
which is impossible according to \cite[Theorem~1.1]{GS1981} since 
\begin{equation*}
    \frac{1}{m-1} = \frac{d}{d-2} < \frac{d+2}{d-2}\,.
\end{equation*}
Consequently, there is $\rho_0\in (0,\infty)$ such that $\mathrm{supp}(P)=\bar{\mathbb{B}}_{\rho_0}$. From~\eqref{6.yd} we infer that
\begin{equation}\label{56}
	P(x)=\left( \frac{m-1}{m} \varphi(x) \right)^{1/(m-1)}\,,\quad x\in \bar{\mathbb{B}}_{\rho_0}\,.
\end{equation}
In addition, the monotonicity and continuity of $\varphi$, along with~\eqref{56}, entail that
\begin{equation*}
	0\le \varphi_+(x) \le \varphi_+\left(\rho_0\frac{x}{|x|}\right) = 0 = \frac{m}{m-1} P(x)\,,\quad x\in \R^d\setminus\bar{\mathbb{B}}_{\rho_0}\,.
\end{equation*}
We infer from the previous two identities that $m P^{m-1}/(m-1) = \varphi_+$ in $\mathbb{R}^d$, hence~\eqref{420a}, which we combine with~\eqref{t4} to obtain~\eqref{varphi3}. Now, invoking~\cite[Theorem~1]{WY2003}, we deduce from~\eqref{b0} and \eqref{varphi3} that $\varphi$ is given by~\eqref{420b}. We have thus established~\eqref{420}.

\medskip

\noindent\textbf{Step~5: $M\ge M_c$.} Finally, the proof that necessarily $M\ge M_c$ is the same as in \cite[Proposition~6.4]{LW26}. Indeed, since $u_l\in\Sigma_a$ for all $l\ge 1$ with $M=a|\mathbb{B}_R|$, we deduce from~\eqref{b2x},~\eqref{420} and~\eqref{N51} that
\begin{align*}
	M & = \lim_{l\to\infty}\int_{\mathbb{B}_{\lambda_l R}}P_l(x)\,\rd x \ge \lim_{l\to\infty}\int_{\mathbb{B}_{\rho_0}} P_l(x)\,\rd x \\ 
	& = \int_{\mathbb{B}_{\rho_0}} P(x)\,\rd x = \frac{1}{\rho_0^{d}} \int_{\mathbb{B}_{\rho_0}} \zeta^{1/(m-1)}\left(\dfrac{x}{\rho_0}\right)\,\rd x = M_c\,,
\end{align*}
which completes the proof of Proposition~\ref{P20}.
\end{proof}

We shall now show that necessarily $M=M_c$, whenever a sequence of radially symmetric stationary states $(u_l,w_l)_{l\ge 1}$ in $\mathcal{S}_{M,\mathrm{rad}}$ drives the Lyapunov functional $\mathcal{L}$ to~$-\infty$. A key difficulty in establishing this result is to obtain sufficiently sharp local estimates to preclude the possibility that mass escapes to infinity along the rescaled sequence~$(P_l)_{l\ge 1}$ introduced and analyzed in Proposition~\ref{P20}. The first step in that direction is a ``sup+inf''-type inequality based on Proposition~\ref{P20}. The arguments are adapted from \cite[Lemma~3]{LS1994} and \cite[Lemma~3]{WY2003}. We omit the proof, as it is the same as that of \cite[Lemma~6.5]{LW26}, to which we refer.

\begin{lem}\label{L25}
Let $a>0$ and 
\begin{equation}
	C_1\ge C_1^* := \dfrac{2(d-2)\sigma_d}{M_c} \left( 1+ \dfrac{m\zeta(0)}{m-1} \right)\,. \label{C1*} 
\end{equation}
Then there is $C_2>0$ such that, for any $u\in\Sigma_a$,
\begin{equation*}
	\phi'(u)(0) + C_1 \inf_{\mathbb{B}_\rho} \phi'(u) \le C_2\rho^{2-d}\,,\quad \rho\in (0,R]\,.
\end{equation*}
\end{lem}

The next step is to substitute the infimum of $\phi'(u)$ over the ball $\mathbb{B}_\rho$ by its value on the boundary $\partial\mathbb{B}_\rho$ for which we require the following result:

\begin{lem}\label{L289}
Let $a>0$ and $u\in\Sigma_a$. There exists $\omega=\omega(d,R)\ge 1$ such that $v=\mathcal{K}[u]$ satisfies
\begin{equation*}
	v(\rho) = \inf_{\partial\mathbb{B}_\rho} v \le \omega \inf_{\mathbb{B}_\rho} v \,,\quad \rho\in (0,R]\,.
\end{equation*}
\end{lem}

\begin{proof}
Set
\begin{equation*}
    g(r) := \Gamma\left(\frac{d}{2}\right) \left(\frac{2}{r}\right)^{(d-2)/2} I_{(d-2)/2}(r), \qquad r>0\,,
\end{equation*}    
where $I_{(d-2)/2}$ is the modified Bessel function of the first kind, and recall that $g$ is the solution to the initial value problem
\begin{equation*}
    -\partial_r \big(r^{d-1}\partial_r g(r)\big)  + r^{d-1} g(r) = 0\,,\quad r>0\,, \qquad g(0)=1\,,\quad g'(0)=0\,,
\end{equation*}
satisfying $g(r)\ge 1$ and $\partial_r g(r)\ge 0$ for $r\ge 0$. Introducing the Wronskian 
\begin{equation*}
    W(r) := r^{d-1}\big(g \partial_r v - v \partial_r g\big)(r)\,, \qquad r\in [0,R]\,, 
\end{equation*}
we infer from the definitions of $v$ and $g$ and the non-negativity of $u$ that
\begin{equation*}
   \partial_r W(r) = r^{d-1} \big(g(v-u) - vg\big)(r) = -r^{d-1} (gu)(r) \le 0, \quad r\in [0,R]\,, 
\end{equation*}
with $W(0)=0$. Hence $W \le 0$ in $[0,R]$, so that $v/g$ is non-increasing on $[0,R]$. Now, consider $\rho\in (0,R]$ and let $r_\rho\in [0,\rho]$ be such that
\begin{equation*}
    v(r_\rho) = \inf_{\mathbb{B}_\rho} v \ge 0\,.
\end{equation*}
Then the monotonicity of $v/g$ and the properties of $g$ imply that
\begin{equation*}
    v(\rho) \le \frac{g(\rho)}{g(r_\rho)} v(r_\rho) \le g(\rho) v(r_\rho) \le g(R) v(r_\rho)\,,
\end{equation*}
from which Lemma~\ref{L289} follows with $\omega:=g(R)$.
\end{proof}

We now establish a local estimate for functions in $\Sigma_a$, adapting an argument designed in \cite[Proposition~4]{WY2003}, see also \cite[Lemma~2]{LS1994}.

\begin{lem}\label{L28}
Let $a>0$ and $M=a|\mathbb{B}_R|$. There is $\delta^*\in (0,1)$ depending only on $d$, $R$, and $M$ with the following property: for any $\delta\in (0,\delta^*]$, there is $C_3(\delta)>0$ depending only on $d$, $R$, $M$,  and $\delta$ such that, for all $u\in \Sigma_a$,
\begin{equation*}
	\phi'(u)(\rho) \le C_3(\delta) \rho^{2-d} - \delta \phi'(u)(0) \,,\quad \rho\in (0,R]\,.
\end{equation*}
\end{lem}

\begin{proof}
Let $\rho\in (0,R]$. Since $v=\mathcal{K}[u] = \phi'(u) - \langle\phi'(u)\rangle+a$, see Lemma~\ref{Lx}, it follows from Lemma~\ref{L289} that 
\begin{align*}
\phi'(u)(\rho) & =  \langle \phi'(u)\rangle -a + v(\rho) = \langle \phi'(u)\rangle -a + \inf_{\partial\mathbb{B}_\rho} v\\
& \le \langle \phi'(u)\rangle -a + \inf_{\mathbb{B}_\rho} v + (\omega-1) \inf_{\mathbb{B}_\rho} v\\
& = \inf_{\mathbb{B}_\rho} \phi'(u) + (\omega-1) \inf_{\mathbb{B}_\rho} v\,.
\end{align*}
For $r\in (0,\rho]$ we infer from~\eqref{n9} that
\begin{align*}
	 v(r) = \phi'(u)(r) -\langle \phi'(u)\rangle +a & \le |\phi'(u)(r) -\langle \phi'(u)\rangle | +a \\
	& \le  \frac{M}{\sigma_d(m-1)}\left(\frac{r^{2-d}}{d}+\frac{R^{2-d}}{2}\right) +a\,,
\end{align*}
so that
\begin{align*}
	\inf_{\mathbb{B}_\rho} v & \le \inf_{0<r<\rho} \left[ \frac{M}{\sigma_d(m-1)} \left( \frac{r^{2-d}}{d}+\frac{R^{2-d}}{2} \right) + a \right]\\
	& \le \frac{M}{\sigma_d(m-1)}\frac{d+2}{d}\rho^{2-d} +a\,.
\end{align*}
Consequently, gathering the outcome of the previous computations yields
\begin{equation*}
	\phi'(u)(\rho) = \inf_{\partial\mathbb{B}_\rho} \phi'(u) \le \inf_{\mathbb{B}_\rho}\phi'(u) + (\omega-1) \left( \frac{M}{\sigma_d}\frac{d+2}{d-2} \rho^{2-d} +a \right)\,.
\end{equation*}
Introducing
\begin{equation*}
	c_3:= (\omega-1) \left( \frac{M}{\sigma_d}\frac{d+2}{d-2} +aR^{d-2}\right)\,,
\end{equation*}
we obtain
\begin{equation} \label{j78}
	\phi'(u)(\rho)\le  \inf_{\mathbb{B}_\rho} \phi'(u)+ c_3 \rho^{2-d}\,,\quad \rho\in (0,R]\,.
\end{equation}
Recalling that $C_1^*$ is defined in~\eqref{C1*}, we set  $\delta^* := 1/C_1^*$ so that $C_1:=1/\delta\ge C_1^*$ for a given $\delta\in (0,\delta^*]$. Therefore, Lemma~\ref{L25} and~\eqref{j78} yield a constant $C_2(\delta)>0$ (depending on $C_1$ and thus on $d$, $R$, $M$, and $\delta$) such that
\begin{align*}
	\phi'(u)(\rho) & \le \delta \left( C_2(\delta) \rho^{2-d} - \phi'(u(0)) \right) + c_3 \rho^{2-d} \\
	& \le \left( \delta C_2(\delta)  + c_3 \right) \rho^{2-d} - \delta \phi'(u)(0)
\end{align*}
for $\rho\in (0,R]$, and the proof is complete with $C_3(\delta) := \delta C_2(\delta)  + c_3$.
\end{proof}

Lemma~\ref{L28} allows us to show that any unbounded sequence in $\Sigma_a$ concentrates at the origin and to determine the associated space scale.

\begin{cor}\label{XX}
Let $a>0$, $M=a|\mathbb{B}_R|$ and consider a sequence $(u_l)_{l\ge 1}$ in $\Sigma_a$ such that $\lambda_l^d:=\|u_l\|_\infty\to \infty$.
Then there exists a constant $\kappa>0$ depending on $d$, $R$, and $M$ such that
\begin{equation}\label{limphi}
	\lim_{l\to\infty} \left(\max_{r\in [\kappa/\lambda_l,R] }u_l(r)\right)=0\,.
\end{equation}
\end{cor}

\begin{proof}
Consider $l\ge 1$ and pick $r_l\in [0,R]$ such that $u_l(r_l) = \lambda_l^d = \|u_l\|_\infty$. Owing to~\eqref{N6}, we find a constant $c_4=c_4(d,R,M)>0$ such that $\lambda_l^d = u_l(r_l) \le c_4 r_l^{-d}$, and thus
\begin{equation}\label{h0}
	r_l\le \frac{c_4^{1/d}}{\lambda_l}\,,\quad l\ge 1\,.
\end{equation}
Introducing $\varrho_l := \lambda_l^{-(d+2)/d}$, we split the analysis according to the ordering of $r_l$ and~$\varrho_l$.

\noindent\textbf{Case~1.} If $r_l\le \varrho_l$, then it follows from~\eqref{n7} that
\begin{align*}
	\lambda_l^d - u_l(0) & = \int_0^{r_l} \partial_r u_l(r)\,\rd r \le \frac{\|u_l\|_\infty}{dm} \int_0^{r_l} r u_l(r)^{2-m}\,\rd r \\
	& \le \frac{\|u_l\|_\infty^{(d+2)/d}}{2dm} r_l^2 \le \lambda_l^{d+2} \varrho_l^2 = \lambda_l^{d-(4/d)}\,.
\end{align*}
\noindent\textbf{Case~2.} If $r_l> \varrho_l$, then 
\begin{equation*}
	\lambda_l^d - u_l(0) = \int_0^{\varrho_l} \partial_r u_l(r)\,\rd r + \int_{\varrho_l}^{r_l} \partial_r u_l(r)\,\rd r\,.
\end{equation*}
On the one hand, it follows from \textbf{Case~1} that
\begin{equation*}
    \int_0^{\varrho_l} \partial_r u_l(r)\,\rd r \le \lambda_l^{d-(4/d)}\,.
\end{equation*}
On the other hand, we use~\eqref{N7A} and H\"older's inequality to obtain
\begin{align*}
    \int_{\varrho_l}^{r_l} \partial_r u_l(r)\,\rd r & \le c_1 \int_{\varrho_l}^{r_l} r^{3-d} u_l(r)^{2-m}\,\rd r = c_1 \int_{\varrho_l}^{r_l} r^{-(d+1)(d-2)/d} \big( r^{d-1} u_l(r)\big)^{2/d}\,\rd r \\
    & \le c_1 \left( \int_{\varrho_l}^{r_l} r^{d-1} u_l(r)\,\rd r \right)^{2/d} \left( \int_{\varrho_l}^{r_l} r^{-(d+1)}\,\rd r \right)^{(d-2)/d} \\
	& \le c_1 \left( \frac{M}{\sigma_d} \right)^{2/d} \left( \frac{\varrho_l^{-d}}{d} \right)^{(d-2)/d} \le c_1 \left( \frac{M}{\sigma_d} \right)^{2/d} \lambda_l^{d-(4/d)}\,.
\end{align*}
Therefore, there is $c_5=c_5(d,R,M)>0$ such that, in both cases, we end up with
\begin{equation*}
    \lambda_l^d - u_l(0) \le c_5 \lambda_l^{d-(4/d)}\,, \qquad l\ge 1\,.
\end{equation*}
Since $\lambda_l\to\infty$,  there is $l_2\ge 1$ such that
\begin{equation*}
	1\le \frac{\lambda_l^d}{2} \le \left( 1-\frac{c_5}{\lambda_l^{4/d}} \right) \lambda_l^d \le u_l(0)\,,\quad l\ge l_2\,.
\end{equation*}
It now follows from~\eqref{L20b}, the identity $d(m-1)=d-2$ and the above lower bound that
\begin{equation*}
	\frac{m}{m-1}\left(2^{1-m}\lambda_l^{d-2}-1\right)\le \phi'(u_l)(0) \,,\quad l\ge l_2\,.
\end{equation*}
As $d\ge 3$ and $\lambda_l\to\infty$, there is $c_6=c_6(d,R,M)>0$ such that, for a possibly larger choice of $l_2$, 
\begin{equation}\label{h3}
	2 c_6 \lambda_l^{d-2}\le \phi'(u_l)(0) \,,\quad l\ge l_2\,.
\end{equation}
Combining~\eqref{h3} with Lemma~\ref{L28} (with $\delta=\delta^*$) ensures that, for $\rho\in (0,R]$ and $l\ge l_2$,
\begin{align*}
	\phi'(u_l)(\rho) & \le C_3(\delta^*) \rho^{2-d} - \delta^* \phi'(u_l)(0) \le C_3(\delta^*) \rho^{2-d} - 2 c_6 \delta^* \lambda_l^{d-2} \\
	& = - c_6 \delta^* \lambda_l^{d-2} + C_3(\delta^*) \lambda_l^{d-2} \left[ (\rho \lambda_l)^{2-d} - \frac{c_6 \delta^*}{C_3(\delta^*)} \right]\,.
\end{align*}
Setting $\kappa^{2-d} := c_6 \delta^*/C_3(\delta^*)$, we infer from the above inequality that there is $l_3\ge l_2$ such that, for $l\ge l_3$, we have $\kappa \lambda_l^{-1} \le R$ and
\begin{equation}\label{h4}
	\phi'(u_l)(\rho) \le - c_6 \delta^* \lambda_l^{d-2}, \quad \rho\in [\kappa/\lambda_l,R]\,.
\end{equation}
In particular, $u_l(\rho)<1$ for $\rho\in [\kappa/\lambda_l,R]$ and $l\ge l_3$ by Lemma~\ref{L20}.

Finally, for $l\ge l_3$ and $\rho\in [\kappa/\lambda_l,R]$,
\begin{align*}
	-\phi'(u_l)(\rho) & = \int_{u_l(\rho)}^1 \phi''(s)\,\rd s = m \int_{u_l(\rho)}^1 \frac{(1+s)^{m-1}}{s}\,\rd s \\
	& \le m 2^{m-1} \big(-\log{(u_l(\rho))}\big)\,.
\end{align*}
Equivalently, thanks to~\eqref{h4}, 
\begin{equation*}
	u_l(\rho) \le \exp\left\{ - \frac{c_6 \delta^*}{m 2^{m-1}} \lambda_l^{d-2} \right\}\,,
\end{equation*}
hence
\begin{equation*}
	\max_{r\in [\kappa/\lambda_l,R/2] } u_l \le \exp\left\{ - \frac{c_6 \delta^*}{m 2^{m-1}} \lambda_l^{d-2} \right\}\,, \quad l\ge l_3\,.
\end{equation*}
As $d\ge 3$ and $\lambda_l\to\infty$, claim~\eqref{limphi} readily follows from the above inequality.
\end{proof}

Corollary~\ref{XX} implies that the Lyapunov functional $\mathcal{L}$ can be unbounded from below along a sequence of radially symmetric stationary solutions in $\mathcal{S}_{M,rad}$ only if $M=M_c$ as stated in the next result.

\begin{prop}\label{P3}
Let $M>0$ and suppose there is a sequence $(u_l,w_l)_{l\ge 1}$ in $\mathcal{S}_{M,rad}$ satisfying
\begin{equation}\label{n39}
	\lim_{l\to\infty} \mathcal{L}(u_l,w_l) = -\infty\,.
\end{equation} 
Then $M=M_c$. In particular, if $0<M\not= M_c$, then
\begin{equation*}
	\mu_M^s := \inf_{(u,w)\in \mathcal{S}_{M,rad}} \mathcal{L}(u,w) > - \infty\,.
\end{equation*}
\end{prop}

\begin{proof}
Proposition~\ref{P3} follows exactly as \cite[Proposition~6.8, Corollary~6.9]{LW26}. 
\end{proof}

Proposition~\ref{P3} is the final ingredient needed to complete the proof of part~\textbf{(II)} of Theorem~\ref{THM1}.

\begin{proof}[Proof of Theorem~\ref{THM1}~\textbf{(II)}]
We argue exactly as in \cite[Theorem~\ref{THM1}~\textbf{(II)}]{LW26}, following the approach initiated in \cite{Ho2002,HW2001}.
\end{proof}

Let us finally mentioned that the equality~\eqref{infl1} is proved by contradiction with the same argument as in the proof of Theorem~\ref{THM1}~\textbf{(II)}, see \cite[Remark~6.10]{LW26}.

\section*{Acknowledgments}

Part of this work was done while PhL enjoyed the hospitality and support of the Institut f\"ur angewandte Mathematik, Leibniz Universit\"at Hannover.

For the purpose of Open Access, a CC-BY license has been applied by the authors to all versions
of this article leading up to the Author Accepted Manuscript.

\bigskip

\paragraph{\textbf{Declaration on the use of AI}} Lemma~\ref{L289} was obtained with the help of Claude Opus 5.5.

\paragraph{\textbf{Data availability}} Our manuscript has no associated data.

\paragraph{\textbf{Conflict of interest}} The authors declare that there is no conflict of interest.

\appendix
\section{Properties of \texorpdfstring{$\phi$}{}}\label{sec.apA}

We recall some properties of the function $\phi: (0,\infty)\rightarrow (0,\infty)$ defined by
\begin{equation}
	\phi''(z)=m\frac{(1+z)^{m-1}}{z}\,,\quad z>0\,,\qquad \phi(1)=\phi'(1)=0\,, \tag{\ref{phi}}
\end{equation}
which were established in \cite[Appendix~A]{LW26}. We begin with upper and lower bounds of $\phi$ in terms of the related power function $s\mapsto s^m$.

\begin{lem}\label{lemA1}
For each $\varepsilon>0$, there is a positive constant $c_\phi(\varepsilon)>0$, depending only on $m$ and $\varepsilon$, such that 
\begin{equation}
    \phi(s) \le \frac{1+\varepsilon}{m-1} s^m + c_\phi(\varepsilon)\,, \quad s\ge 0\,. \label{eq.A1}
\end{equation}
In addition,
\begin{equation}
    \phi(s)\ge \frac{s^m}{m-1} + 1 - \frac{m}{m-1} s\,, \quad s\ge 0\,. \label{eq.A2}
\end{equation}
\end{lem}

We next recall estimates involving $\phi$ and $\phi'$, which we use throughout the analysis of stationary solutions to~\eqref{E}.

\begin{lem}\label{L20}
The function $\phi$ obeys the inequalities
\begin{equation}\label{L20a}
	s\phi'(s)\le m\phi(s) +\phi''(1)(s-1)\,,\quad s> 0\,,
\end{equation}
and
\begin{equation}\label{L20b}
	\frac{m}{m-1}(s^{m-1}-1)\le \phi'(s)\le \frac{m 2^{m-1}}{m-1} (s^{m-1}-1)\,,\quad s\ge 1\,,
\end{equation}
with $\phi'(s)<0$ for $s\in (0,1)$.
\end{lem}

We omit the proofs of Lemmas~\ref{lemA1} and~\ref{L20} and refer to \cite[Appendix~A]{LW26}. We end up this section with an additional property of $\phi$, allowing us to study the boundedness from below of the Lyapunov function in Section~\ref{sec.4}.

\begin{lem}\label{LTJ}
The function $r\mapsto 2\phi(r)-r\phi'(r)$ is bounded from below on $(0,\infty)$.
\end{lem}

\begin{proof}
First, $\phi'(r) \le 0$ for $r\in (0,1]$ so that $2\phi(r)- r\phi'(r)\ge 2\phi(r)\ge 0$ for $r\in (0,1]$. Next, it follows from~\eqref{eq.A2} and~\eqref{L20a} that, for $r\ge 1$,
\begin{align*}
    2\phi(r)-r\phi'(r) & = (2-m) \phi(r) + m \phi(r) - r\phi'(r) \ge (2-m) \phi(r) -\phi''(1)(r-1)\\
    & \ge \frac{2-m}{m-1} r^m + 2-m - \frac{m(2-m)}{m-1} r -\phi''(1)(r-1),
\end{align*}
which is clearly bounded from below in view of $m>1$.
\end{proof}

\bibliographystyle{siam}
\bibliography{Lit_Chem}

\end{document}